\documentclass[a4paper,11pt]{article}

\usepackage[latin1]{inputenc} 
\usepackage[T1]{fontenc} 
\usepackage{lmodern} 
\usepackage[left=3.2cm,top=3.5cm,right=3.2cm,nohead,foot=1cm]{geometry}
\usepackage[english]{babel} 
\usepackage{tikz} 
\usetikzlibrary{arrows,shadows}

\usepackage{amsfonts,latexsym,amssymb} 
\usepackage{mathrsfs,amsmath}
\usepackage[all]{xy}

\usepackage{xcolor} 

\usepackage[numbers]{natbib} 

\newlength{\defbaselineskip}
\newcommand{\Proof}{{\noindent {\bf Proof:   }}}

\newcommand{\EndProof}{{\hfill $\Box \qquad $ \endtrivlist}\par \vspace{0.5cm}}%

\renewcommand{\phi}{\varphi}

\newcommand{\vecA}{\mathbf{A}}

\newcommand{\vecE}{\mathbf{E}}
\newcommand{\vecF}{\mathbf{F}}
\newcommand{\vecG}{\mathbf{G}}
\newcommand{\vecI}{\mathbf{I}}
\newcommand{\vecK}{\mathbf{K}}

\newcommand{\vecO}{\mathbf{O}}

\newcommand{\vecR}{\mathbf{R}}
\newcommand{\vecS}{\mathbf{S}}
\newcommand{\vecT}{\mathbf{T}}
\newcommand{\vecU}{\mathbf{U}}

\newcommand{\vecW}{\mathbf{W}}

\newcommand{\vecId}{\mathbf{Id}}
\newcommand{\vecc}{\mathbf{c}}
\newcommand{\vece}{\mathbf{e}}
\newcommand{\vecf}{\mathbf{f}}

\newcommand{\vecn}{\mathbf{n}}
\newcommand{\vecp}{\mathbf{p}}

\newcommand{\vecr}{\mathbf{r}}
\newcommand{\vecs}{\mathbf{s}}
\newcommand{\vecu}{\mathbf{u}}
\newcommand{\vecv}{\mathbf{v}}
\newcommand{\vecw}{\mathbf{w}}
\newcommand{\vecx}{\mathbf{x}}
\newcommand{\vecy}{\mathbf{y}}

\newcommand{\vecsigma}{\boldsymbol\sigma}
\newcommand{\vecxi}{\boldsymbol\xi}

\newcommand{\Div}{\mbox{ div }}
\newcommand{\ds}{\displaystyle}

\newcommand{\B}{ \mathtt{B} }

\newenvironment{Remarque} {  \normalfont \begin{flushright} $\rightsquigarrow$ \vline \- \begin{minipage}{10cm} \textbf{brouillon.}}
                          {  \end{minipage} \end{flushright}  }

\newcommand{\comment}[1]{}

    \newtheorem{theorem}{Theorem}[section]
    \newtheorem{lemma}[theorem]{Lemma}
    \newtheorem{corollary}[theorem]{Corollary}
     
    \newtheorem{proposition}[theorem]{Proposition}

    \newtheorem{remark}{Remark}[section]

\renewcommand{\comment}[1]{ {#1} }
\renewcommand{\comment}[1]{ {#1} }

\title{Enhanced controllability of low Reynolds number swimmers in the presence of a wall
}

\author{François Alouges\footnote{Centre de Mathématiques Appliquées de l'Ecole Polytechnique, CMAP, Ecole Polytechnique, Plateau de Saclay, 91128 Palaiseau, France, Tel.: +33 1 69 33 46 31, francois.alouges@polytechnique.edu}, Laetitia Giraldi\footnote{Centre de Mathématiques Appliquées de l'Ecole Polytechnique, CMAP, Ecole Polytechnique, Plateau de Saclay, 91128 Palaiseau, France,
              Tel.: +33 1 69 33 46 09, laetitia.giraldi@polytechnique.edu\vspace{0.1cm}}}

\renewcommand{\vec}[1]{\overrightarrow{#1}}

\begin{document}
\maketitle

		\newcount\hh
		\newcount\mm
		\mm=\time
		\hh=\time
		\divide\hh by 60
		\divide\mm by 60
		\multiply\mm by 60
		\mm=-\mm
		\advance\mm by \time
		\def\hhmm{\number\hh:\ifnum\mm<10{}0\fi\number\mm}

\textbf{Abstract} :
Swimming, i.e., being able to advance in the absence of external forces by performing cyclic shape changes, is particularly demanding at low Reynolds numbers which is the regime of interest for micro-organisms and micro-robots. We focus on self-propelled stokesian robots composed of assemblies of balls and we prove that the presence of a wall has an effect on their motility. To rest on what has been done in \cite{AlougesDeSimone10} for such systems swimming on $\mathbf{R}^3$, we demonstrate that a controllable swimmer remains controllable in a half space whereas the reachable set of a non fully controllable one is affected by the presence of a wall.
\\



\section{Introduction}

Swimming at low Reynolds number is now a well established topic of research which probably dates back to the pioneering work of Taylor \cite{Taylor51} who explains how a micro-organism can swim without inertia. Later on, Purcell \cite{Purcell77} formalized the so-called ``scallop theorem'' which states that, due to the reversibility of the viscous flow, a reciprocal deformation of the body cannot lead to a displacement of the swimmer. However, this obstruction can be circumvented using many swimming strategies \cite{Purcell77}. Swimmers can be distinguished with respect to their ability to change their shape or to impose rotational motions of some parts of their body in order to create viscous friction forces on the fluid, and produce by reaction, the propulsion.

Many applications  are concerned by this problem as for example, the conception of medical micro devices. The book by J.P. Sauvage \cite{Sauvage01} presents a lot of engine models adapted for tiny devices while
the design and fabrication of such engines have been recently investigated by e.g. B. Watson, J. Friend, and L. Yeo \cite{WatsonFriend09}. As an example, let us quote the toroidal swimmer, first introduced by Purcell \cite{Purcell77} and which has been subsequently improved by A.M Leshansky and O. Kenneth \cite{LeshanskyKenneth08}, Y. Or and M. Murray \cite{OrMurray09}, A. Najafi and R. Zargar \cite{NajafiZargar10} among others.
 
The strategy for swimming consists in a cyclic deformation of body with a non-reciprocal motion. The first swimmer prototype belonging to this class is the three link swimmer also designed by Purcell \cite{Purcell77}. More recently, R. Golestanian and A. Ajdari \cite{GolestanianAjdari08} introduced the Three-sphere swimmer which is geometically simpler and allows for exact calculations of motion and speed \cite{AlougesDeSimone08}, or even explicit in some asymptotic regimes \cite{GolestanianAjdari08}.

In the continuation of \cite{AlougesDeSimone08}, F. Alouges, A. DeSimone, L. Heltai, A. Lefebvre, and B. Merlet  \cite{AlougesDeSimone10}  showed that the trajectory of the Three-sphere swimmer is governed by a differential equation whose control functions correspond to the rate of changing shape. The swimming capability of the device now is recast in terms of a control problem to which classical results apply.

Of particular importance for applications is the issue of the influence of any boundary on the effective swimming capabilities of micro-devices or real micro-swimmers. Indeed, boundaries clearly affect the hydrodynamics and may have an influence on the swimmer's capabilities. In that direction, an biological study of Rothschild \cite{Rothschild63} claimed for instance that spermatozoids tend to accumulate on walls. More recently,  H. Winet, G. S. Bernstein, and J. Head \cite{WinetBernstein84} proved this related boundary effect for the sperm of humans which evolves in a narrow channel. Swimming in a geometrically confined environment then became a subject of major interest, in particular to model this attraction phenomenon (see \cite{SmithBlake10}, \cite{GaffneyGadelha11},\cite{BerkeH.-C.-Berg08}). D. J Smith, E. A. Gaffney and J. R. Blake \cite{SmithGaffney09} have described the motion of a stylized bacterium propelled by a single flagellum and they show that the attraction by the wall is effective. Later, H. Shum, E.A Gaffney, and J. Smith  \cite{ShumGaffney10} investigated to which extent this attraction effect is impacted by a change in swimmer's morphology. 


On a more theoretical side, other approaches provide results that show an attraction effect by the wall. A. P. Berke and P. Allison \cite{BerkeAllison08}, modelling the swimmer with a simple dipole, put in evidence an attraction due to the presence of the wall. Y. Or and M. Murray \cite{OrMurray09} derived the swimmer dynamics near a wall for three various swimmers, but with unvarying shapes. 
The case of a changing shape swimmer has been studied by R. Zargar and A. Najafi \cite{ZargarNajafi09}, where the dynamics of the Three-sphere swimmer in the presence of a wall is given. However, some fundamental symmetry are not satisfied in their swimmer's motion equation. 
 
The aim of this paper is to attack the same problem (the influence of a plane wall in the motion of the swimmer) by means of control theory. 
Several recent works present a controllability results for a self-propelled micro-swimmers in a space without boundary, as example, let us quote the paper of J. Lohéac, J. F. Scheid and M. Tucsnak \cite{LoheacScheid11} and the study of J. Lohéac and A. Munnier \cite{LoheacMunnier12} made of the spherical swimmer in the whole space (see also \cite{Chambrion11} for the same kind of results in a perfect fluid). 
Furthermore, F. Alouges, A. DeSimone, L. Heltai, A. Lefebvre, and B. Merlet  \cite{AlougesDeSimone10} deal with the controllability on $\mathbf{R}^3$ for the Three-sphere swimmer and others specific swimmers. The question that we want to address now is whether the presence of the plane wall modifies the controllability results. We here prove two results in that direction. Namely, considering the fully controllable Four sphere swimmer proposed in \cite{AlougesDeSimone10}, we show that in the half space, the swimmer remains fully controllable, while a Three-sphere swimmer enriches its reachable set, at least generically which seems at first sight contradictory with earlier results. Indeed, although previous works show an attraction from the boundary, the set of reachable points could be of higher dimension. In other words, if the dynamics is somehow more constrained due to the presence of the wall, the set of points that the swimmer may reach could be larger than what it was without the wall.

%


The rest of this paper is organized as follow. In Section \ref{Swimmers}, we describe the two model swimmers to which our analytical and numerical tools are later applied. Section~\ref{MainResults} presents the main controllability results associated with the introduced swimmers. In Section \ref{MathsSetting}, we show  that swimming is indeed an affine control problem without drift by using a similar approach than F. Alouges, A. DeSimone, L. Heltai, A. Lefebvre, and B. Merlet in \cite{AlougesDeSimone10}. 
The controllability result  is proved in Section \ref{FourSphere} for the Four-sphere swimmer and in Section \ref{ThreeSphere}  for the Three-sphere swimmer. Concluding remarks are given in Section \ref{Conclusion}.

%
\section{Swimmers}
\label{Swimmers}

We carry on the study of specific swimmers that were considered in  \cite{AlougesDeSimone10} in $\mathbf{R}^3$. In order to fix notation, the wall is modeled by the plane $W=\{(x,y,z)\in \mathbb{R}^3\mbox{ s. t. } y = 0\}$, and the swimmers, which consist of $N$ spheres $(B_i)_{i=1..N}$ of radii $a$ connected by thin jacks, are assumed to move in the half space $\mathbf{R}^3_+=\{(x,y,z)\in \mathbf{R}^3\mbox{ s. t. } y> 0\}$. As in \cite{AlougesDeSimone10} , the viscous resistance associated with the jacks is neglected and the fluid is thus assumed to fill the whole set $\mathbf{R}^3_+\setminus \cup_{i=1}^N B_i$. The state of the swimmer is described by two sets of variables :
\begin{itemize}
\item the shape variables, denoted by $\vecxi$ (here in $\mathbf{R}^{N-1}$ or $\mathbf{R}^{N}$), which define the lengths of the jacks. A stroke consists in changing the lengths of these jacks in a periodic manner ;
\item the position variables, denoted by $\vecp \in \mathbf{R}_+^3 \times SO(3)$, which define swimmer's position and orientation in the half-space.
\end{itemize}

In what follows, we call $\mathcal{S} \subset \mathbf{R}^M$ for a suitable $M \in \mathbf{N}$ the set of admissible states $(\vecxi,\vecp)$ that we assume to be a connected nonempty 
smooth submanifold of $\mathbf{R}^M$.  We thereafter focus on two swimmers that have been considered in the literature, the Three-sphere swimmer (see \cite{NajafiGolestanian04}, \cite{AlougesDeSimone08}, \cite{AlougesDeSimone10}) and the Four sphere swimmer (see \cite{AlougesDeSimone10}). It turns out that this latter is easier to understand that the former, and we therefore start with it.

\subsection{The Four-sphere swimmer}
\label{4Sswimmer}
\begin{figure}[ht]
\centering
\includegraphics[scale=1]{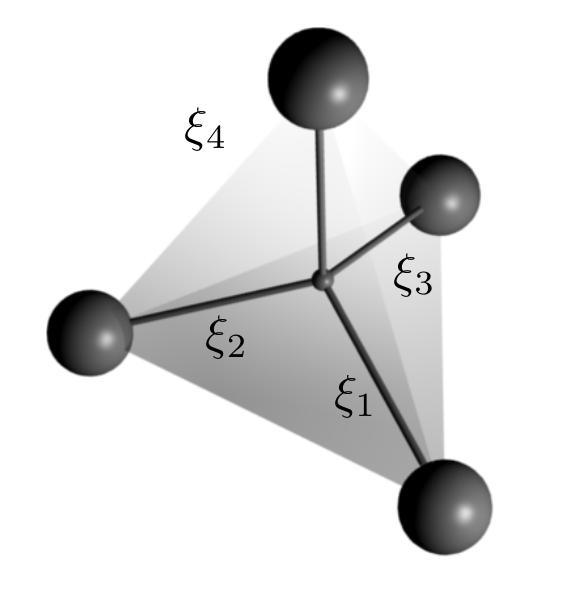}
\caption{The Four-sphere swimmer.\label{FourSphereSwimmer}}
\end{figure}

We consider a regular tetrahedron $(\vecS_1,\vecS_2,\vecS_3,\vecS_4)$ with center $\vecO \in \mathbf{R}^3_+$. 
The swimmer consists on four balls linked by four arms of fixed directions $\vec{\vecO\vecS_i}$ which are able to elongate and shrink (in a referential associated to the swimmer). 
The four ball cluster is completely described by the list of parameters $(\vecxi,\vecp) = (\xi_1, \dots, \xi_4,\vecc,\boldsymbol{\mathcal{R}}) \in \mathcal{S}= (\sqrt{\frac{3}{2}}a,\infty)^4 \times \mathbf{R}^3_+ \times SO(3)  \subset (\sqrt{\frac{3}{2}}a,\infty)^4\times \mathbf{R}^{6}$ (a rotation $\boldsymbol{\mathcal{R}}\in SO(3)$ is uniquely characterized by its 3 dimensional rotaion vector). It is known (see  \cite{AlougesDeSimone10}) that the Four sphere swimmer is controllable in $\mathbf{R}^3$.
This means that it is able to move to any point and with any orientation under the constraint of being self-propelled, and when the surrounding flow is dominated by the viscosity. This swimmer is depicted in Fig. \ref{FourSphereSwimmer}.

\newpage
\subsection{The Three-sphere swimmer}
\label{ThreesphereIntroduction}
\begin{figure}[ht]
\centering
\begin{center}
\begin{tikzpicture}[scale=0.7]
\draw (-2,-2) circle(1);
\draw (0,0) circle (1) ;
\draw (0,0)--(4,4);
\draw (4,4) circle(1);
\draw (-2,-2)--(0,0);
\draw (-2,-2) node[below]{$\vecx_1$};
\draw (0,0) node[below right]{$\vecx_2$};
\draw (4,4)  node[below right]{$\vecx_3$};
\draw [>=latex,->] (0,0) -- (135 : 1);
\draw (0.2,0.5) node[left]{$a$};
\draw (-0.1,0.1) -- (0.1,-0.1);
\draw (-1.9,-2.1) -- (-2.1,-1.9);
\draw (3.9,4.1) -- (4.1,3.9);
\draw [dashed][->,>=stealth] (0,0) -- (5,0);
\draw (2,0) arc (0:45:2) ;
\draw (22.5:2) node[right]{$\theta$} ;
\draw [<->](-1,1) -- (3,5);
\draw [<->](-1,1) -- (-3,-1);
\draw (-2.25,0.3) node[]{$\xi_1$};
\draw (0.5,3) node[]{$\xi_2$};
\end{tikzpicture}
\end{center}
\caption{The Three-sphere swimmer.\label{ThreeSphereSwimmer}}
\end{figure}
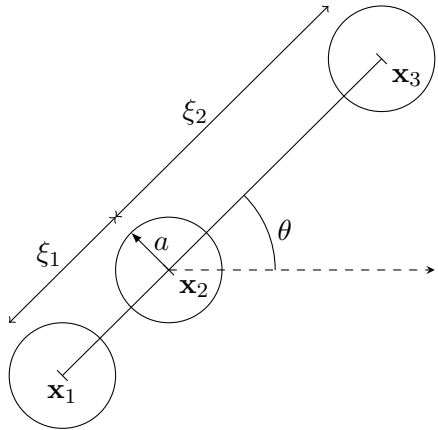

This swimmer is composed of three aligned spheres as shown in Fig. \ref{ThreeSphereSwimmer}.
We assume that at $t=0$ the swimmer starts in the vertical half-plane  $H= \{(x,y,0)\in \mathbf{R}^3\mbox{ s. t. }z=0,y\geq 0\}$, it is clear from the symmetry of the problem that the swimmer stays in $H$ for all time, for whatever deformation of its arms it may carry out. We characterize swimmer's position and orientation in $H$ by the coordinates $(\vecc,\theta) \in  \mathbf{R}^2 \times [0,2\pi]$, where $\vecc \in H$ is the position of one of the three spheres, and $\theta$ is the angle between the swimmer and the $x-$axis. Therefore, in that case, the swimmer is completely described by the vector $(\vecxi,\vecp) = (\xi_1,\xi_2,\vecc,\theta)\in \mathcal{S}=(2a,\infty)^2 \times H \times [0,2\pi) \subset (2a,\infty)^2 \times \mathbf{R}^2\times \mathbf{R}/{2\pi\mathbf{Z}}$. 
In the three dimensional space $\mathbf{R}^3$ (when there is no boundary), it is obvious by symmetry that the angle $\theta$ cannot change in time, and thus this swimmer is not fully controllable. One of the main contributions of this paper is to understand the modifications of this behavior due to the presence of the plane wall.

\section{The main results}
\label{MainResults}


Consider any of the swimmers described in the previous Sections, and assume it is self-propelled in a three dimensional half space viscous flow modeled by Stokes equations.
In this paper, we will establish that both swimmers are locally fully controllable almost everywhere in $\mathcal{S}$. By this we mean the precise following statements.

\begin{theorem}
\label{Thm1}
Consider the Four-sphere swimmer described in Section \ref{ThreesphereIntroduction}, and assume it is self-propelled in a three dimensional viscous flow modeled by Stokes equations in the half space $\mathbf{R}^3_+$. Then for almost any initial configuration $(\vecxi^i,\vecp^i)\in \mathcal{S}$, any final configuration $(\vecxi^f,\vecp^f)$ in a suitable neighborhood of $(\vecxi^i,\vecp^i)$ and any final time $T>0$, there exists a stroke $\vecxi\in \mathcal{W}^{1,\infty}([0,T])$, satisfying
$\vecxi(0)=\vecxi^i$ and $\vecxi(T)=\vecxi^f$ and such that if the self-propelled swimmer starts in position $\vecp^i$ with the shape $\vecxi^i$ at time $t=0$, it ends at position $\vecp^f$ and shape $\vecxi^f$ at time $t=T$ by changing its shape along $\vecxi(t)$.
\end{theorem}

\begin{theorem}
\label{Thm2}
Consider the Three-sphere swimmer described in Section \ref{4Sswimmer}, and assume it is self-propelled in a three dimensional viscous flow modeled by Stokes equations in the half space $\mathbf{R}^3_+$. Then for almost any initial configuration $(\vecxi^i,\vecp^i)\in \mathcal{S}$ such that $\vecp^i \in H$, any final configuration $(\vecxi^f,\vecp^f)$ in a suitable neighborhood of $(\vecxi^i,\vecp^i)$ with $\vecp^f \in H$ and any final time $T>0$, there exists a stroke $\vecxi\in \mathcal{W}^{1,\infty}([0,T])$, satisfying
$\vecxi(0)=\vecxi^i$ and $\vecxi(T)=\vecxi^f$ and such that if the self-propelled swimmer starts in position $\vecp^i$ with the shape $\vecxi^i$ at time $t=0$, it ends at position $\vecp^f$ and shape $\vecxi^f$ at time $t=T$ by changing its shape along $\vecxi(t)$ and staying in $H$ for all time $t\in [0,T]$.

\end{theorem} 

\begin{remark}
The sense of ``almost every initial configuration'' can be further precised as everywhere outside a (possibly empty) analytic manifold of codimension 1.
\label{rmkanal}
\end{remark}

%
%
%

The proof of the controllability of the Four-sphere swimmer is given in Section $\ref{FourSphere}$ whereas Section $\ref{ThreeSphere}$ is devoted to demonstrate devoted to demonstrate Theorem ~\ref{Thm2}. 


\section{Mathematical setting of the problem}
\label{MathsSetting}

As for their 3D counterparts, the equation of motion of both swimmers take the form of an affine control problem without drift. In this section, we detail the derivation of this system.

\subsection{Modelization of the fluid}

The flow takes place at low Reynolds number and we assume that inertia of both the swimmer and the fluid is negligible. As a consequence, denoting by $\Omega= \cup_{i=1}^N B_{i}$ the space occupied by the swimmer, the flow in $\mathbf{R}^3_+\setminus \Omega$ satisfies the (static) Stokes equation

\begin{equation}
\label{Stokes}
\left\{
\begin{array}{lllll}
- \mu \Delta \vecu + \nabla p = 0 \,\,\, \textrm{in $\mathbf{R}^3_+ \setminus \Omega$,}\\
\Div \,\vecu = 0 \,\,\, \textrm{in $\mathbf{R}^3_+ \setminus \Omega$,}\\
- \vecsigma \vecn = \vecf   \,\,\, \textrm{on $\partial \Omega$,}\\
\vecu = 0 \,\,\, \textrm{on $\partial \mathbf{R}^3_+$,}\\
\vecu \to 0  \,\,\, \textrm{at $\infty$.}\\
\end{array}
\right.
\end{equation}

Here, we have denoted by $\vecsigma = \mu ( \nabla \vecu + \nabla^t \vecu) - p \mathbf{Id}$ the Cauchy stress tensor, $\vecn$ is the unit normal to $\partial \Omega$ pointing outward to the swimmer. We also set
$$
\mathcal{V}=\{ \vecu \in \mathcal{D}' (\mathbf{R}^3_+\setminus \Omega,\mathbf{R}^3)\, \vert \, \nabla \vecu \in L^2(\mathbf{R}^3_+\setminus \Omega), \, \frac{\vecu(\vecr)}{\sqrt{1+\vert \vecr \vert^2}} \in L^2(\mathbf{R}^3_+\setminus \Omega) \}\,.
$$
It is well known that $\mathcal{V}$ is a Hilbert space when endowed with the norm (and the associated scalar product)
$$
\Vert \vecu \Vert_{\mathcal{V}}^2 := \int_{\mathbf{R}^3_+\setminus \Omega} \vert \nabla \vecu \vert^2\,.
$$ 
\newline

We also assume that $\vecf \in H^{-1/2}(\partial \Omega)$ in order to obtain a unique solution $(\vecu,p)$ to the problem \eqref{Stokes} in $\mathcal{V} \times L^2(\mathbb{R}^3_+\setminus \Omega)$ which can be expressed in terms of the associated Green's function (obtained by the method of ``images'', see \cite{Blake71}) as 

\begin{equation}
\label{RepresentationIntegrale}
\vecu(\vecr) = \int_{\partial\Omega} \vecK(\vecr, \mathbf{s}) \vecf(\mathbf{s}) d\mathbf{s},
\end{equation}
where the matricial Green function $\vecK = (K_{ij})_{i,j=1,2,3}$ is given by
\begin{equation}
\label{DefK}
\vecK(\vecr,\vecr_0)=\vecG(\vecr-\vecr_0)+\vecK_1(\vecr,\vecr_0)+\vecK_2(\vecr,\vecr_0)+\vecK_3(\vecr,\vecr_0)\,,
\end{equation}
the four functions $\vecG$, $\vecK_1$, $\vecK_2$ and $\vecK_3$ being respectively
the Stokeslet
\begin{equation}
\label{stokeslet}
\vecG(\vecr) = \frac{1}{8\pi \mu} \left( \frac{\vecId }{| \vecr  |} + \frac{\vecr \otimes \vecr}{\vert \vecr \vert^3} \right) 
\end{equation}
and the three ``ìmages''
\begin{equation}
\label{DefK1}
\vecK_1(\vecr,\vecr_0) = -\ds \frac{1}{8 \pi \mu} \left( \frac{\vecId}{\vert \vecr' \vert} + \frac{\vecr'\otimes \vecr'}{\vert \vecr' \vert^3} \right)\,, 
\end{equation}
\begin{equation}
\label{DefK2}
K_{2,ij}(\vecr,\vecr_0) = \frac{1}{4\pi \mu}  y_0^2 \left(1 - \delta_{j2}\right) \left( \frac{\delta_{ij}}{\vert \vecr' \vert^3} - \frac{3 r'_i r'_j}{\vert \vecr' \vert^5} \right)\,,
\end{equation}
\begin{equation}
\label{DefK3}
K_{3,ij}(\vecr,\vecr_0) = -\frac{1}{4\pi \mu}  y_0 \left(1 - 2 \delta_{j2} \right) \left( \frac{r'_2}{\vert \vecr' \vert^3}\delta_{ij} - \frac{r'_j}{\vert \vecr' \vert^3}\delta_{i2} + \frac{r'_i}{\vert \vecr' \vert^3} \delta_{j2} - \frac{3 r'_i r'_j r'_2}{\vert \vecr' \vert^5} \right) \,. 
\end{equation}
Here $\vecr_0=(x_0,y_0,z_0)$ and $\vecr' = \vecr - \tilde{\vecr}_0$, where $\tilde{\vecr}_0=(x_0,-y_0,z_0)$ stands for the ``image'' of $\vecr_0$, that is to say, the point symmetric to $\vecr_0$ with respect to the wall. 


Let $B$ be the sphere of radius $1$ centered at the origin.
We identify the boundary of the domain occupied by the swimmer, $\partial \Omega$, with $(\partial B)^N$ and we represent by $\vecf_i \in H^{-1/2}( \partial B)$ the distribution of force on the sphere $B_i$. Correspondingly, $\vecu_i \in H^{1/2}( \partial B)$ stands for the velocity distribution on the sphere $B_i$ (and of the fluid due to non-slip contact).

Following \cite{AlougesDeSimone10}, we denote by $\mathcal{T}_{(\vecxi,\vecp)}$ the Neumann-to-Dirichlet map

\begin{equation}
\label{NeumannDirichlet}
\begin{array}{cccc}
\mathcal{T}_{(\vecxi,\vecp)} : & \mathcal{H}^{-1/2} & \to & \mathcal{H}^{1/2}\\
& (\vecf_1, \dots, \vecf_N) & \mapsto& (\vecu_1, \dots, \vecu_N)\\
\end{array}
\end{equation}
where we have denoted by $\mathcal{H}^{\pm1/2}$ the space $(H^{\pm1/2}( \partial B))^N$. It is well known that the map $\mathcal{T}_{(\vecxi,\vecp)}$ is a one to one mapping onto while its inverse is continuous. 

Using \eqref{RepresentationIntegrale}, we can express $\vecu_i$ ($i=1,2,3$) by
\begin{equation}
\label{VitesseBoule}
\begin{array}{cccc}
\forall \vecr \in \partial B, & \vecu_i(\vecr) & = &\ds \sum_{j=1}^{N} \ds \int_{\partial B} \vecK(\vecx_i +a \vecr,\vecx_j + a\vecs) \vecf_j(\vecs) d\vecs \\
 & &:= & \ds \sum_{j=1}^{N} \langle \vecf_j, \vecK(\vecx_i + a\vecr, \vecx_j  + a\cdot )  \rangle_{\partial B}\,,\\  
\end{array}   
\end{equation} 
where $\langle \cdot, \cdot \rangle_{\partial B}$ stands for the duality $\big(H^{-1/2}( \partial B ),H^{1/2}( \partial B)\big)$.

\begin{proposition}
\label{AnalycityProp}
The mapping $(\vecxi,\vecp) \mapsto \mathcal{T}_{(\vecxi,\vecp)}$ is analytic from $\mathcal{S}$ into $ \mathcal{L}(\mathcal{H}^{-1/2},\mathcal{H}^{1/2})$.
Furthermore, $\mathcal{T}_{(\vecxi,\vecp)}$ is an isomorphism for every $(\vecxi,\vecp) \in \mathcal{S}$, and the mapping $(\vecxi,\vecp) \mapsto \mathcal{T}^{-1}_{(\vecxi,\vecp)}$ is also analytic.
\end{proposition}

\Proof
The proof is identical to the one given in \cite{AlougesDeSimone10}, replacing the the Stokeslet by the Green kernel $\vecK$ which is also analytic outside its singularity.
\EndProof

\begin{remark}
As the direct consequence, the mapping $\mathcal{T}_{(\vecxi,\vecp)}$ and its inverse depends analytically on $a$.
\end{remark}

\subsection{Equation of motion}

In this section, we use the self-propulsion assumption in order to express the dynamics of the swimmer as an affine control system without drift.

\begin{proposition}
There exists a family of vectorfields $\vecF_i \in \mathbf{T}\mathcal{S}$, such that the state of the swimmer is described by the following ODE, 
\begin{equation}
\label{EqMotion}
\frac{d}{dt} \left( \begin{array}{ll} \vecxi \\ \vecp \end{array} \right) = \sum_i\vecF_i(\vecxi,\vecp) \dot{\xi}_i\,.
\end{equation}

\end{proposition}

\Proof

This equation of motion is by now classical in this context (see \cite{AlougesDeSimone10}, \cite{AlougesDeSimone08}, \cite{Dal-MasoDeSimone10} or \cite{Munnier09}).
Let us recall the principle of its derivation.

At any time $t$, the swimmer occupies a domain $\Omega_t$ (we therefore denote  by $\Omega_0$ the domain occupied by the swimmer at time $t\,=\,0$).
 We also define the map $\Phi$ which associates to the points of $\partial \Omega_0 \times \mathcal{S}$, the current point in $\partial \Omega_t$,

\begin{equation}
\begin{array}{cccc}
\Phi : & \partial \Omega_0 \times \mathcal{S} & \to & \partial \Omega_t \\
 & (\vecx_0, \vecxi, \vecp) & \mapsto & \vecx_t \,.
\end{array}
\end{equation}
When inertia is negligible, self-propulsion of the swimmer implies that the total viscous force and torque exerted by the surrounding fluid on the swimmer vanish i.e.,
\begin{equation}
\label{ForceTorque}
\left\{
\begin{array}{ccccc}
\vecF:=& \displaystyle \int_{\partial \Omega_t} \mathcal{T}^{-1}_{\vecp,\vecxi}\left( \frac{\partial \Phi}{\partial t} \right) d\vecx_t & = &\displaystyle \int_{\partial \Omega_t} \mathcal{T}^{-1}_{\vecp,\vecxi}\left((\partial_{\vecp} \Phi) \dot{\vecp}+(\partial_{\vecxi \Phi}) \dot{\vecxi}\right) d\vecx_t &=0\,,\\
\vecT:=& \displaystyle \int_{\partial \Omega_t} \vecx_t \times \mathcal{T}^{-1}_{\vecp,\vecxi}\left( \frac{\partial \Phi}{\partial t}\right) d\vecx_t & = &\displaystyle \int_{\partial \Omega_t} \vecx_t \times \mathcal{T}^{-1}_{\vecp,\vecxi}\left((\partial_{\vecp} \Phi) \dot{\vecp}+(\partial_{\vecxi} \Phi) \dot{\vecxi}\right) d\vecx_t &=0\,.
\end{array}
\right.
\end{equation}

From the linearity of the Neumann-to-Dirichlet map, we deduce that the system (\ref{ForceTorque}) reads as a linear system which depends on $\dot{\vecp}$ and $\dot{\vecxi}$. By inverting it, we get $\dot{\vecp}$ linearly in terms of $\dot{\vecxi}$. Assuming that $\vecxi \in \mathbf{R}^k$ for some $k \in \mathbf{N}$, we thus obtain

\begin{equation}
\dot{\vecp}= \displaystyle \sum_{i=1}^k \vecW_i(\vecxi,\vecp) \dot{\xi}_i\,,
\end{equation}
which becomes (\ref{EqMotion}) when we call $\vecF_i :=  \left( \begin{array}{c}
 \vecE_i \\ \vecW_i
\end{array}
\right)$, where $\vecE_i$ is the $i$-th vector of the canonical basis. 
\EndProof

Let us recall some notations which are used to study the controllability of such systems of ODE (see for instance \cite{Jurdjevic97}).

 Let $F$ and $G$ be two vector fields defined on a smooth finite dimensional manifold $\mathcal{M}$. The Lie bracket of $F$ and $G$ is the vector field defined at any point $X\in \mathcal{M}$ by $[F,G](X):=\left( F \cdot \nabla \right) G(X) - \left( G \cdot \nabla \right) F(X)$. For a family of vector fields  $\mathcal{F}$ on $\mathcal{M}$, $Lie(\mathcal{F})$ denotes the Lie algebra generated by $\mathcal{F}$. Namely, this is the smallest algebra - defined by the Lie bracket operation - which contains $\mathcal{F}$  (therefore $\mathcal{F}\subset Lie(\mathcal{F})$ and for any two vectorfields $F \in Lie(\mathcal{F})$ and $G\in Lie(\mathcal{F})$, the Lie bracket $[F,G] \in Lie(\mathcal{F})$). 
Eventually, for any point $X \in \mathcal{M}$, $Lie_X(\mathcal{F})$ denotes the set of all tangent vectors $V(X)$ with $V$ in $Lie(\mathcal{F})$. It follows that $Lie_X(\mathcal{F})$ is a linear subspace of $T_X \mathcal{M}$ and is hence finite-dimensional. 

Lie brackets and Lie algebras play a prominent role in finite dimensional control theory. Indeed, we recall Chow's theorem:
\begin{theorem}[Chow \cite{Coron56}]
\label{ChowGlobal}
Let $\mathcal{M}$ be a connected nonempty manifold. 
Let us assume that $\mathcal{F}= \left(F_i\right)_{i=1}^m$, a family of vector fields on $\mathcal{M}$, is such that  $F_i \in C^{\infty}(\mathcal{M},T\mathcal{M}) \,, \forall i \in\{1,\cdots,m\}\,.$\\
Let us also assume that $$Lie_X(\mathcal{F} )= T_X(\mathcal{M})  \,,\, \forall X \in \mathcal{M}\,.$$\\
Then, for every $(X^0,X^1) \in \mathcal{M} \times \mathcal{M}$, and for every $T >0$, there exists $u \in L^{\infty}([0,T];\mathbf{R}^m)$ such that the solution of the Cauchy problem,

\begin{equation}
\left\{
\begin{array}{ll}
\dot{X} = \displaystyle \sum_{i=1}^m u_i F_i(X)\,,\\
X(0) = X^0\,,
\end{array}
\right.
\end{equation}
is defined on $[0,T]$ and satisfies $X(T) = X^1$.
\end{theorem}

The theorem \ref{ChowGlobal} is a global controllability result, we also recall the one which gives a small-time local controllability. 

\begin{theorem}[\cite{Coron56}, p. 135]
Let $\Omega$ be an nonempty open subset of $\mathbf{R}^n$, that $\mathcal{F}= \left(F_i\right)_{i=1}^m$, a family of vector fields, such that $F_i \in C^{\infty}(\Omega,\mathbf{R}^n) \,, \forall i \in\{1,\cdots,m\}\,.$\\
Let $X_{e}$ such that $$Lie_{X_e}(\mathcal{F} )= \mathbf{R}^n \,.$$\\
Then, for every $\epsilon >0$, there exists a real number $\eta>0$ such that, for every $(X^0,X^1) \in \left\{X \,s.\,\,t.\,\|X-X_e\|<\eta\right\}^2$, there exists a bounded measurable function $u: \left[0,\epsilon\right] \to \mathbf{R}^n$ such that the solution of the Cauchy problem
\begin{equation}
\left\{
\begin{array}{ll}
\dot{X} = \displaystyle \sum_{i=1}^m u_i F_i(X)\,,\\
X(0) = X^0\,,
\end{array}
\right.
\end{equation}
 is defined on $[0,\epsilon]$ and satisfies $X(\epsilon) = X^1$.
\end{theorem}

When the vector fields are furthermore analytic (and the manifold $\mathcal{M}$ is also analytic) one also has the Hermann-Nagano Theorem of which we will make an important use in the theoretical study of the controllability for our model swimmers.
\begin{theorem}[Hermann-Nagano \cite{Jurdjevic97}]
Let $\mathcal{M}$ be an analytic manifold, and $\mathcal{F}$ a family of analytic vector fields on $\mathcal{M}$. Then
\begin{enumerate}
\item each orbit of $\mathcal{F}$ is an analytic submanifold of $\mathcal{M}$, and
\item if $\mathcal{N}$ is an orbit of $\mathcal{F}$, then the tangent space of $\mathcal{N}$ at $X$ is given by $Lie_X(\mathcal{F})$. In particular, the dimension of $Lie_X(\mathcal{F})$ is constant as $X$ varies over $\mathcal{N}$.
\end{enumerate}
\end{theorem}

In our context, the family of vector fields is given by $\mathcal{F}=(\vecF_i)_{1\leq i\leq k}$  which are defined on the manifold $\mathcal{M}=\mathcal{S}$, and
the controls $u_i$ are given by the rate of shape changes $\dot{\xi}_i$. 
In view of the preceding theorems, the controllability question of our model swimmers raised by Theorems \ref{Thm1} and \ref{Thm2} 
relies on the dimension of the Lie algebra
generated by the vectorfields $(\vecF_i)_{1\leq i\leq k}$ which define the dynamics of the swimmer. In particular they are direct consequences of the following Lemma.

\begin{lemma}
For almost every point (in the sense of remark \ref{rmkanal}) $(\vecxi,\vecp)\in \mathcal{S}$, the Lie algebra generated by the vectorfields $(\vecF_i)_{1\leq i\leq k}$ at $(\vecxi,\vecp)$ is equal to $T_{(\vecxi,\vecp)}\mathcal{S}$. 
\label{lemmalie}
\end{lemma}

The proof of this lemma is developed until the rest 
of the paper. Several tools are used in order to characterize this dimension among which we mainly use asymptotic behavior and
symbolic computations. As we shall see, although the theory is clear, the explicit computation (or at least asymptotic expressions) is by no means
obvious and requires a lot of care. In particular, before using symbolic calculations, a rigorous proof of the expansion, together with a 
careful control of the remainders in the expressions allowed us to go further.

\section{The Four-sphere swimmer}
\label{FourSphere}

In this section, we give the proof of the controllability result stated in Theorem~\ref{Thm1}. 
\Proof
The argument of the proof is based on the fact that $\vecK$ given by (\ref{DefK}) satisfies
\begin{equation}
\vecK(\vecr,\vecr') = \vecG(\vecr-\vecr') + O\left(\frac{1}{y}\right)\,,
\end{equation}
where $\vecr=(x,y,z)$ and $\vecr'=(x',y',z')$ are two points of $\mathbf{R}^3_+$, and $\vecG$ is the Green function of the Stokes problem in the whole space $\mathbf{R}^3$, namely the Stokeslet,  defined by (\ref{stokeslet}).



As a consequence, we obtain that the Neumann to Dirichlet map given by (\ref{NeumannDirichlet}) satisfies for a swimmer of shape $\vecxi$ at position $\vecp=(p_x,p_y,p_z,\mathbf{\mathcal{R}}) \in \mathbf{R}^3_+ \times SO_3$
\begin{equation}
\mathcal{T}_{(\vecxi,\vecp)} = \mathcal{T}^0_{\vecxi} + O\left(\frac{1}{p_y}\right)\,,
\end{equation} 
where $\mathcal{T}^0_{\vecxi}$ is the Neumann-to-Dirichlet map associated to the Green function $\vecG\,$.

The system (\ref{ForceTorque}) now reads 
\begin{equation}
\left\{
\begin{array}{ll}
\displaystyle \int_{\partial \Omega_t} \left(  \left(\mathcal{T}_{\vecxi}^0\right)^{-1} + \,O\left(\frac{1}{p_y}\right) \right)\left((\partial_\vecp \Phi) \dot{\vecp}+ (\partial_{\vecxi}  \Phi) \dot{\vecxi}\right) d\vecx_t &  =0\,,\\
\displaystyle \int_{\partial \Omega_t} \vecx_t \times \left( \left(\mathcal{T}_{\vecxi}^0\right)^{-1} + O\left(\frac{1}{p_y}\right)\right) \left((\partial_{\vecp} \Phi) \dot{\vecp}+(\partial_{\vecxi} \Phi) \dot{\vecxi}\right) d\vecx_t& \, =0\,.
\end{array}
\right.
\end{equation}
Consequently, the ODE (\ref{EqMotion}) becomes
%
\begin{equation}
\label{EqMotion2}
\frac{d}{dt} \left( \begin{array}{ll} \vecxi \\ \vecp \end{array} \right) = \sum_{i=1}^4 \left(\vecF^0_i(\vecxi)+O\left(\frac{1}{p_y}\right)\right) \dot{\xi}_i\,,
\end{equation}
where $(\vecF^0_i)_{i=1,\cdots,4}$ are the vector fields obtained in the case of the whole space $\mathbf{R}^3$.
 
In other words, we obtain the convergence 
\begin{equation}
\label{ConvergenceSimple1}
 \vecF(\vecxi,\vecp) =  \vecF^0(\vecxi) + O\left(\frac{1}{p_y}\right)\mbox{ as }p_y\rightarrow +\infty
\end{equation}
and also for all its derivatives to any order.
%

It has been proved in \cite{AlougesDeSimone10} that $\mbox{dim} \,\,Lie_{\vecxi}(\vecF^0) = 10$ at all admissible shape $\vecxi$, showing the global controllability in the whole space of the underlying swimmer. We thus obtain that for $p_y$ sufficiently large 
\begin{equation}
\label{ConditionDimension}
\mbox{dim} \,\,Lie_{(\vecxi,\vecp)}(\vecF) = 10\,,
\end{equation}
and therefore due to the analyticity of the vector fields $(\vecF_i)_{i=1,\cdots,4}$, (\ref{ConditionDimension}) holds in a dense subset of $\mathcal{S}$. This shows that the system satisfies the full rank condition almost everywhere in $\mathcal{S}$ and proves Lemma \ref{lemmalie} in this context, and thus Theorem \ref{Thm1} by a simple application of Chow's theorem.
\EndProof

The preceding proof can be generalized to any swimmer for which the Lie algebra satisfies the full rank condition in $\mathbf{R}^3$. We now turn to an example for which this is not the case, namely the Three-sphere swimmer of Najafi Golestanian \cite{NajafiGolestanian04}. Indeed, when there is no boundary, this swimmer is constrained to move along its axis of symmetry. The purpose of the next section is to understand to which extent this is still the case when there is a flat boundary.

%
%
%
%
%
%
%
\section{The Three-sphere swimmer}
\label{ThreeSphere}

This section details the proof of Theorem ~\ref{Thm2}. It is organized in several subsections, each of them focusing on a particular step of the proof.
In the subsection \ref{Subsec:EquMotion}, by introducing some notations, we recall the expression of the equation of motion of the Three-sphere swimmer. Subsection  \ref{Subsec:Symmetry} deals with the special symmetry which have to be verify by the vector fields of the motion equation. From this symmetry properties, we deduce the reachable set of the particular case where the swimmer is perpendicular to the wall.  In the subsection \ref{Subsec:Approx}, we give an expansion of the Neumann-To-Dirichlet mapping associated to the Three-sphere swimmer and its inverse, in the case where the radius of the sphere $a$ is small enough and the distance of the arm is sufficiently large. In subsection \ref{Subsec:Self-propulsion},  we deduce from this previous approximation an expansion of the motion equation of the swimmer, for $a$ sufficiently small. Finally, the subsection \ref{Subsec:Expansion} presents some formal calculations of the vectors fields of the motion equation and their Lie brackets which leads to obtain, almost everywhere, the dimension of its Lie algebra.  

\subsection{Equation of motion for the Three-sphere swimmer}
 \label{Subsec:EquMotion}

From Section $\ref{ThreesphereIntroduction}$, we know that the swimmer's position is parameterized by the vector $(x,y,\theta)$ where $(x,y)$ is the coordinate of the center of $B_2$ as depicted in Fig. \ref{ThreeSphereSwimmer} and $\theta$ is the angle between the swimmer and the $x-$axis. We recall that $\vecxi:=(\xi_1,\xi_2)$  stands for the lengths of both arms of the swimmer.

The motion equation  (\ref{EqMotion}) thus reads,
\begin{equation}
\label{EqMotionTreeSphere}
\frac{d}{dt} \left( \begin{array}{ll} \xi_1 \\ \xi_2 \\ x \\ y \\ \theta \end{array} \right) = \vecF_1(\vecxi,x,y,\theta) \dot{\xi}_1 + \vecF_2(\vecxi,x,y,\theta) \dot{\xi}_2\,.
\end{equation}
Notice that, from translational invariance of the problem, both $\vecF_1$ and $\vecF_2$ actually do not depend on $x$.

In what follows, we denote by 
$$
d_{(\vecxi,y,\theta)}=\mbox{dim }Lie_{(\vecxi,y,\theta)}(\vecF_1,\vecF_2)
$$ 
the dimension of the Lie algebra $Lie_{(\vecxi,y,\theta)}(\vecF_1,\vecF_2)\subset  \mathbb{R}^5$ at $(\vecxi,y,\theta)$. It is clear, since $\vecF_1$ and $\vecF_2$ are independent one to another and never vanish, that
\begin{equation}
\label{EncadrementDimension}
2\leqslant d_{(\vecxi,y,\theta)}\leqslant 5\,.
\end{equation}

\subsection{Symmetry properties of the vector fields}
 \label{Subsec:Symmetry}

\begin{proposition}
\label{Symmetry1}
Let $\vecS$ be the $5\times 5$ matrix defined by
$$
\vecS=\left(
\begin{array}{ccccc}
0 & 1 & 0 & 0 & 0\\
1 & 0 & 0 & 0 & 0\\
0 & 0 & -1 & 0 & 0\\
0 & 0 & 0 & 1 & 0\\
0 & 0 & 0 & 0 & -1
\end{array}
\right)\,.
$$
Then one has for all $\vecxi = (\xi_1,\xi_2,x,y,\theta)\in \mathcal{S}$ 
\begin{equation}
\vecF_1(\xi_1,\xi_2,y,\theta) = \vecS \vecF_2(\xi_2,\xi_1,y,2\pi-\theta)
\end{equation}
and similarly for the Lie bracket 
\begin{equation}
[\vecF_1,\vecF_2](\xi_1,\xi_2,y,\theta) = \vecS [\vecF_2,\vecF_1](\xi_2,\xi_1,y,2\pi-\theta)\,.
\end{equation}
\end{proposition}

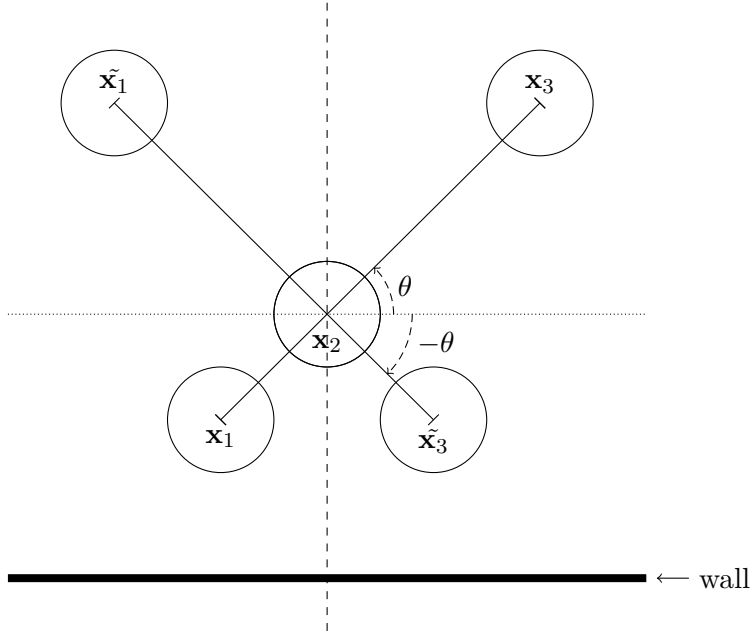
\begin{figure}[ht]
\centering
\begin{center}
\begin{tikzpicture}[scale=0.7]
\draw (-2,-2) circle(1);
\draw (0,0) circle (1) ;
\draw (0,0)--(4,4);
\draw (4,4) circle(1);
\draw (-2,-2)--(0,0);
\draw (-2,-2) node[below]{$\vecx_1$};
\draw (0,-0.2) node[below]{$\vecx_2$};
\draw (4,4)  node[above]{$\vecx_3$};
\draw (-1.9,-2.1) -- (-2.1,-1.9);
\draw (3.9,4.1) -- (4.1,3.9);
\draw (2,-2) circle(1);
\draw (0,0) circle (1) ;
\draw (0,0)--(-4,4);
\draw (-4,4) circle(1);
\draw (2,-2)--(0,0);
\draw (-1.9,-2.1) -- (-2.1,-1.9);
\draw (3.9,4.1) -- (4.1,3.9);
\draw (-4.1,3.9) --(-3.9,4.1);
\draw (1.9,-2.1) --(2.1,-1.9);
\draw (2,-2) node[below]{$\tilde{\vecx_3}$};
\draw (-4,4) node[above]{$\tilde{\vecx_1}$};
\draw [dashed] (0,-6)--(0,6);
\draw [line width=3pt] (-6,-5)--(6,-5);
\draw [<-] (6.2,-5)--(6.8,-5);
\draw (6.8,-5) node[right]{$\textrm{wall}$};
\draw [densely dotted] (-6,0)--(6,0);
\draw [densely dashed,->] (1.25,0) arc (0:45:1.25); 
\draw [densely dashed,->] (1.6,0) arc (0:-45:1.6); 
\draw (25:1.25) node[right]{$\theta$};
\draw (-20:1.6) node[right]{$-\theta$};
\end{tikzpicture}
\end{center}
\caption{The plane symmetry which links the situation at $(\xi_1,\xi_2,x,y,\theta)$ with those at $(\xi_2,\xi_1,x,y,2\pi-\theta)$. In both cases, solutions to Stokes flow are also symmetric one to another. \label{FigSymmetry}
}
\end{figure}

\Proof
Although the plane breaks the 3D axisymmetry along the swimmer's axis, we can still make use of the symmetry with respect to the vertical plane that passes through the center of the first sphere $B_2$. A swimmer with position $(x,y,\theta)$ and shape $(\xi_1,\xi_2)$ is transformed to one at position $(x,y,2\pi-\theta)$ and shape $(\xi_2,\xi_1)$ (see Fig. \ref{FigSymmetry}). Making use of the fact that corresponding solutions to Stokes equations are symmetric one to another, we easily get the proposition.


Eventually, one deduces the Lie bracket symmetries by applying the former symmetries on the vectorfields themselves. An easy recurrence shows that the same identities hold for any Lie bracket of any order of the vectorfields $\vecF_1$ and $\vecF_2$. In particular one has for instance
\begin{equation}
[\vecF_1,[\vecF_1,\vecF_2]](\xi_1,\xi_2,y,\theta) = \vecS [\vecF_2,[\vecF_2,\vecF_1]](\xi_2,\xi_1,y,2\pi-\theta)\,.
\end{equation}

\EndProof

As a direct consequence of proposition $\ref{Symmetry1}$, we deduce that the fourth coordinate of the Lie bracket $[\vecF_1,\vecF_2]$ vanishes at $(\xi,\xi,y,0)$ and at $(\xi,\xi,y,\pi)$.
  
\begin{proposition}
\label{Symmetry2}
Let $\vecT$ be the $5\times 5$ matrix defined by
$$
\vecT=\left(
\begin{array}{ccccc}
1 & 0 & 0 & 0 & 0\\
0 & 1 & 0 & 0 & 0\\
0 & 0 & -1 & 0 & 0\\
0 & 0 & 0 & 1 & 0\\
0 & 0 & 0 & 0 &  -1
\end{array}
\right)\,.
$$
Then one has for all $\vecxi = (\xi_1,\xi_2,x,y,\theta)\in \mathcal{S}$, and $i =1,2$
\begin{equation}
\vecF_i(\xi_1,\xi_2,y,\theta) = \vecT \vecF_i(\xi_1,\xi_2,y,\pi-\theta)
\end{equation}
and similarly for the Lie bracket 
\begin{equation}
[\vecF_1,\vecF_2](\xi_1,\xi_2,y,\theta) = \vecT [\vecF_1,\vecF_2](\xi_1,\xi_2,y,\pi-\theta)\,.
\end{equation}

%
\end{proposition}

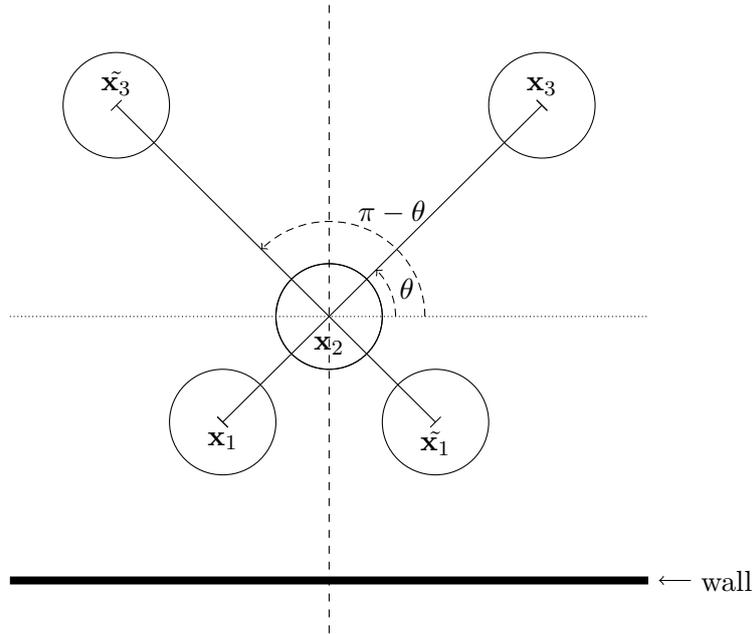
\begin{figure}[ht]
\centering
\begin{center}
\begin{tikzpicture}[scale=0.7]
\draw (-2,-2) circle(1);
\draw (0,0) circle (1) ;
\draw (0,0)--(4,4);
\draw (4,4) circle(1);
\draw (-2,-2)--(0,0);
\draw (-2,-2) node[below]{$\vecx_1$};
\draw (0,-0.2) node[below]{$\vecx_2$};
\draw (4,4)  node[above]{$\vecx_3$};
\draw (-1.9,-2.1) -- (-2.1,-1.9);
\draw (3.9,4.1) -- (4.1,3.9);
\draw (2,-2) circle(1);
\draw (0,0) circle (1) ;
\draw (0,0)--(-4,4);
\draw (-4,4) circle(1);
\draw (2,-2)--(0,0);
\draw (-1.9,-2.1) -- (-2.1,-1.9);
\draw (3.9,4.1) -- (4.1,3.9);
\draw (-4.1,3.9) --(-3.9,4.1);
\draw (1.9,-2.1) --(2.1,-1.9);
\draw (2,-2) node[below]{$\tilde{\vecx_1}$};
\draw (-4,4) node[above]{$\tilde{\vecx_3}$};
\draw [dashed] (0,-6)--(0,6);
\draw [line width=3pt] (-6,-5)--(6,-5);
\draw [<-] (6.2,-5)--(6.8,-5);
\draw (6.8,-5) node[right]{$\textrm{wall}$};
\draw [densely dotted] (-6,0)--(6,0);
\draw [densely dashed,->] (1.25,0) arc (0:45:1.25); 
\draw [densely dashed,->] (1.8,0) arc (0:135:1.8); 
\draw (25:1.25) node[right]{$\theta$};
\draw (80:2) node[right]{$\pi-\theta$};
\end{tikzpicture}
\end{center}
\caption{The plane symmetry which links the situation at $(\xi_1,\xi_2,x,y,\theta)$ with those at $(\xi_1,\xi_2,x,y,\pi-\theta)$. In both cases, solutions to Stokes flow are also symmetric one to another. \label{FigSymmetry2}
}
\end{figure}

\Proof
The two identities readily come from the symmetry which transforms a swimmer with position $(x,y,\theta)$ and a shape $(\xi_1,\xi_2)$ to one at position $(x,y,\theta)$ with the same shape (see Fig. \ref{FigSymmetry2}).

Eventually, one deduces the Lie bracket symmetries by applying the former symmetries on the vectorfields themselves. An easy recurrence shows that the same identities hold for any Lie bracket of any order of the vectorfields $\vecF_1$ and $\vecF_2$.
\EndProof

As a result, in the case where $\theta=\pm\frac{\pi}{2}$, we get the dimension of the Lie algebra of the vector field $\vecF_1$ and $\vecF_2$. 
\begin{corollary}
\label{Corollpi/2}
The dimension of the Lie algebra $Lie_{(\xi_1,\xi_2,y,\pi/2)}(\vecF_1,\vecF_2)$ is less than or equal to $3$.
\end{corollary}

\Proof
We deduce from the preceding proposition that for $i=1,2$ and $j=3,5$, $\vecF^j_i(\xi_1,\xi_2,y,\pm\pi/2)=0$. This simply means that a swimmer starting in the vertical position cannot change its angle $\theta$ and its abscissa $x$ by changing the size of its arms. As a matter of fact, the same holds true for any Lie bracket of $\vecF_1$ and $\vecF_2$ at any order, and we can deduce from this that
$$
d_{(\xi_1,\xi_2,y,\pi/2)}\leq 3\,,
$$
 since any vector of the Lie algebra $Lie_{(\xi_1,\xi_2,y,\pi/2)}(\vecF_1,\vecF_2)$ has a vanishing third and fifth component.
 
Moreover, by using the argument introduced in Section $\ref{FourSphere}$, we get, for almost every $y$, that the dimension of the Lie algebra is almost equal to the one without boundary (see \cite{AlougesDeSimone10}).
\EndProof

\begin{remark}
We denote, for all $x \in \mathbf{R}$, the set 
$$\mathcal{N}_{x} := \left\{(\xi_1,\xi_2,x_0,y,\pm \frac{\pi}{2}) \,\,s.\,\,t. \,\,(\xi_1,\xi_2) \in (2a,\infty)^2 \, y>0 \right\}
$$ 
which corresponds to the case where the swimmer is perpendicular to the wall and $\mathcal{A}_x$ the set of states where the  the dimension of the Lie algebra generated by  $\vecF_1$ and $\vecF_2$ is equal to two, i.e.,
$$
\mathcal{A}_x:=\left\{((\xi_1,\xi_2,x,y,\pm \frac{\pi}{2}) \,s.\,\,t.\, d_{(\xi_1,\xi_2,y,\pi/2)}=2\right\}\,.$$ 
By using the property of analyticity \ref{AnalycityProp}, $\mathcal{A}_x$ is a finite union,
$$\mathcal{A}_x =\ds \bigcup_{y \in F_y} \left\{(\xi_1,\xi_2,x,y,\pm \frac{\pi}{2})\,s.\,\,t.\,(\xi_1,\xi_2) \in (2a,\infty)^2\right\}\,.
$$ 
We deduce, that for all $x$, $\mathcal{N}_{x}\backslash\mathcal{A}_x$ defines a set of three-dimensional orbits strictly included in the manifold $\mathcal{S}$. 
\end{remark}

Furthermore, the proof of the corollary \ref{Corollpi/2} can be applied to all generic positions then, it implies that the dimension of the Lie algebra is almost equal to $3$, almost everywhere, i.e.,
$$
3\leqslant d_{(\xi_1,\xi_2,y,\theta)}\leqslant 5\,.
$$

\subsection{Approximation for small spheres and large distances}
 \label{Subsec:Approx}
For the general case ($\theta \ne \pi/2$), the preceding computation is not sufficient to conclude. In order to proceed, we make an expansion of the
vectorfields and their Lie brackets with respect to $a$ (the radius of the balls) near 0.
 
This part is devoted to the proof of the expansion of the Neumann to Dirichlet map (\ref{expansion}) together with its inverse (\ref{inverse}) at large arms' lengths. Let us first define for all $(i,j)\in \{1,2,3\}^2$, the linear map $\mathcal{T}_{i,j}$ as

\begin{eqnarray*}
\mathcal{T}_{i,j} : H^{-1/2}(\partial B) & \rightarrow &H^{1/2}(\partial B)\\
\vecf_j &\mapsto & \int_{\partial B} \vecK(\vecx_i+ a\cdot,\vecx_j + a\vecs) \, \vecf_j (\vecs) \,\mbox{d}\vecs\,.
\end{eqnarray*}
We recall that the Green kernel $\vecK$ writes (following (\ref{DefK})) as
$$
\vecK(\vecr,\vecr')=\vecG(\vecr-\vecr')+\vecK_1(\vecr,\vecr')+\vecK_2(\vecr,\vecr')+\vecK_3(\vecr,\vecr')\,,
$$
where $\vecG$ is the Stokeslet (see (\ref{stokeslet})) and each kernel is given by the corresponding counterpart in (\ref{DefK}). Eventually, we call ${\mathcal{T}}^0$ the Neumann to Dirichlet map associated to $\vecG$
\begin{eqnarray*}
{\mathcal{T}}^0 : H^{-1/2}(\partial B) & \rightarrow &H^{1/2}(\partial B)\\
\vecf &\mapsto & \int_{\partial B} \vecG( a(\cdot - \vecs)) \, \vecf (\vecs) \,\mbox{d}\vecs\,.
\end{eqnarray*}

\begin{proposition}
\label{propexp}
Let $(i,j)\in \{1,2,3\}^2$. We have the following expansions, valid for $a \ll 1$:
\begin{itemize}
\item if $i\ne j$ then
\begin{equation}
\mathcal{T}_{i,j} = \vecK(\vecx_i,\vecx_j)   \langle \vecf_j, \vecId \rangle_{\partial B} +\vecR_1
\label{idifj}
\end{equation}
where $\ds ||\vecR_1||_{\mathcal{L}(H^{-1/2},H^{1/2})}=O\left(a\right)$\,,
\item otherwise
\begin{equation}
\mathcal{T}_{i,i} = {\mathcal{T}}^0+\sum_{k=1}^3 \vecK_k(\vecx_i,\vecx_i)   \langle \vecf_i, \vecId \rangle_{\partial B} +\vecR_2
\label{ieqj}
\end{equation}
where $\ds ||\vecR_2||_{\mathcal{L}(H^{-1/2},H^{1/2})}=O\left(a\right)$\,. 
\end{itemize}
\end{proposition}

\Proof
Let $(i,j)\in \{1,2,3\}^2$ be such that $i\ne j$, and $\vecf_j \in H^{-1/2}(\partial B)$. We define
\begin{equation}
\forall \vecr\in \partial B\,,\,\,\vecu_i(\vecr):= (\mathcal{T}_{i,j} \vecf_j)(\vecr) = \int_{\partial B} \vecK(\vecx_i+a\vecr,\vecx_j+a\vecs)\vecf_j(\vecs)\mbox{d}\vecs\,,
\label{Tij}
\end{equation}
and
$$
\vecv_i(\vecr)=\vecu_i(\vecr)-\vecK(\vecx_i,\vecx_j)\int_{\partial B}\vecf_j(\vecs)\,d\vecs =\int_{\partial B} \left(\vecK(\vecx_i+a\vecr,\vecx_j+a\vecs)-\vecK(\vecx_i,\vecx_j)\right)\vecf_j(\vecs)\mbox{d}\vecs\,.
$$ 
Our aim is to estimate the $H^{1/2}(\partial B)$ norm of $\vecv_i$. But\footnote{Here and in the sequel, we use the definition for the $H^{1/2}(\partial B)$ norm:
$$
\|\vecv\|_{H^{1/2}(\partial B)}=\min_{\vecw \in H^1(\B,\mathbb{R}^3),\,\vecw=\vecv \mbox{ on }\partial B}\|\vecw\|_{H^1(B)}\,.
$$}
$$
\|\vecv_i\|_{H^{1/2}(\partial B)} \leq \|\vecv_i\|_{H^1(B)},
$$ 
and since $\vecK(\vecx,\vecy)$ is a smooth function in the neighborhood of $\vecx=\vecx_i$ and $\vecy=\vecx_j$, one has $\forall \vecr,\vecs \in B$
\begin{equation}
\left|\vecK(\vecx_i+a\vecr,\vecx_j+a\vecs) - \vecK(\vecx_i,\vecx_j)\right| = O\left(a\right)\,,
\end{equation}
and for the gradients in both $\vecr$ and $\vecs$
\begin{eqnarray*}
&&\left|\nabla_{\vecr} \vecK(\vecx_i+a\vecr,\vecx_j+a\vecs) \right| = O\left(a\right)\,,\\
&&\left|\nabla_{\vecs} \vecK(\vecx_i+a\vecr,\vecx_j+a\vecs) \right| = O\left(a\right)\,,\\
&&\left|\nabla_{\vecr}\nabla_{\vecs} \vecK(\vecx_i+a\vecr,\vecx_j+a\vecs) \right| = O\left(a^2\right)\,.
\end{eqnarray*}
Therefore, we obtain $\forall \vecr \in B$
\begin{eqnarray*}
|\vecv_i(\vecr)|&\leq& \|\vecK(\vecx_i+a\vecr,\vecx_j+a\cdot)-\vecK(\vecx_i,\vecx_j)\|_{H^\frac12}\|\vecf_j\|_{H^{-\frac12}}\\
&\leq& O\left(a\right)\|\vecf_j\|_{H^{-\frac12}}\,,
\end{eqnarray*}
and similarly
\begin{eqnarray*}
|\nabla_{\vecr}\vecv_i(\vecr)|&\leq& \|\nabla_{\vecr} \left(\vecK(\vecx_i+a\vecr,\vecx_j+a\cdot)\right)\|_{H^\frac12}\|\vecf_j\|_{H^{-\frac12}}\\
&\leq& O\left(a\right)\|\vecf_j\|_{H^{-\frac12}}\,.
\end{eqnarray*}
This enables us to estimate the $H^\frac12$ norm of $\vecv_i$ on $\partial B$
\begin{eqnarray*}
\left\|\vecv_i\right\|_{H^\frac12(B)} &\leq& \left\|\vecv_i\right\|_{H^1(B)}\\
 &=& \left(\left\|\vecv_i\right\|^2_{L^2(B)}+\left\|\nabla \vecv_i\right\|^2_{L^2(B)}\right)^\frac12\\
 &\leq&O\left(a\right)\|\vecf_j\|_{H^{-\frac12}}\,. 
\end{eqnarray*}
which proves (\ref{idifj}).

In order to prove (\ref{ieqj}), we use the decomposition (\ref{DefK}) where none of the kernels $(\vecK_i)_{i=1,2,3}$ is singular. Therefore $\forall \vecr\in \partial B$
\begin{eqnarray*}
\vecu_i(\vecr):=( \mathcal{T}_{i,i} \vecf_i)(\vecr) &=& \int_{\partial B} \vecK(\vecx_i+a\vecr,\vecx_i+a\vecs)\vecf_i(\vecs)\mbox{d}\vecs\,\\
&=& \int_{\partial B} \vecG(a(\vecr-\vecs))\vecf_i(\vecs)\mbox{d}\vecs+\int_{\partial B} (\vecK_1+\vecK_2+\vecK_3)(\vecx_i+a\vecr,\vecx_i+a\vecs)\vecf_i(\vecs)\mbox{d}\vecs\\
&=&{\mathcal{T}}^0 \vecf_i + \int_{\partial B} (\vecK_1+\vecK_2+\vecK_3)(\vecx_i+a\vecr,\vecx_i+a\vecs)\vecf_i(\vecs)\mbox{d}\vecs\,.
\end{eqnarray*}
We finish as before, having remarked that for $l=1,2,3$
\begin{equation}
\vecK_l(\vecx_i+a\vecr,\vecx_i+a\vecs) = \vecK_l(\vecx_i,\vecx_i) + O\left(a\right)\,.
\end{equation}
\EndProof


\begin{proposition}
For every $\vecf \in \mathcal{H}^{-1/2}$,
\begin{equation}
\left(\mathcal{T}_{\vecx} \vecf \right)_i (\vecr) = {\mathcal{T}}^0 \vecf_i + \sum_{l=1}^3 \vecK_l(\vecx_i,\vecx_i)   \langle \vecf_i, \vecId \rangle_{\partial B}  + \displaystyle \sum_{j\neq i} \vecK(\vecx_i,\vecx_j) \langle \vecf_j, \vecId\rangle_{\partial B} + \mathcal{R}_i(\vecf),
\label{expansion}
\end{equation}
with $\ds \Vert \mathcal{R}_i\Vert_{\mathcal{L}(\mathcal{H}^{-1/2},\mathcal{H}^{1/2})} =O\left(a\right)$.
\end{proposition}

\Proof

For all $i \in {1,2,3}$, and all $\vecr\in \partial B$
\begin{eqnarray*}
\left(\mathcal{T}_{\vecx} \vecf \right)_i (\vecr)&:=& \displaystyle \int_{\partial B} \vecK(\vecx_i +a\vecr, \vecx_i + a\vecs) \, \vecf_i(\vecs) \mbox{ds} + \sum_{i \neq j} \int_{\partial B} \vecK(\vecx_i +a\vecr, \vecx_j + a\vecs) \, \vecf_j(\vecs) \mbox{d}\vecs \\
&=& \mathcal{T}_{i,i}\vecf_i + \sum_{j \neq i} \mathcal{T}_{i,j} \vecf_j
\end{eqnarray*}
and the result follows from the application of (\ref{idifj}) and (\ref{ieqj}) of Proposition \ref{propexp}.

\EndProof

\begin{proposition}
In the regime $\ds a \ll 1$, one has for every $\vecu \in \mathcal{H}^{1/2}$,
\begin{equation}
\label{inverse}
\begin{array}{ll}
\left(\mathcal{T}^{-1}_{\vecx} \vecu \right)_i = & ({{\mathcal{T}}^0})^{-1} \left( \vecu_i - \sum_{k=1}^3  \vecK_k(\vecx_i,\vecx_i) \langle ({{\mathcal{T}}^0})^{-1} \vecu_i, \vecId \rangle_{\partial B} \right) -\\
& \hspace*{2cm} \displaystyle ({{\mathcal{T}}^0})^{-1} \left(  \sum_{j\neq i} \vecK(\vecx_i,\vecx_j) \langle ({{\mathcal{T}}^0})^{-1} \vecu_j, \vecId \rangle_{\partial B}\right) + \tilde{\mathcal{R}}_i(\vecu) 
\end{array}
\end{equation}
with $\ds \Vert \tilde{\mathcal{R}}_i\Vert_{\mathcal{L}(\mathcal{H}^{1/2},\mathcal{H}^{-1/2})} =O\left(a^3\right)$.
\end{proposition}

\Proof
We recall that
\begin{eqnarray*}
{\mathcal{T}}^0 : H^{-\frac12}(\partial B) &\rightarrow& H^{\frac12}(\partial B)\\
\vecf &\mapsto& \int_{\partial B} \vecG(a(\cdot-\vecs)) \vecf(\vecs)\,d\vecs\,,
\end{eqnarray*}
and define for $l=1,2,3$ the operators
\begin{eqnarray*}
\mathcal{S}_l : H^{-\frac12}(\partial B) &\rightarrow& H^{\frac12}(\partial B)\\
\vecf &\mapsto& \int_{\partial B} \vecK_l(\vecx_i,\vecx_i) \vecf(\vecs)\,d\vecs\,,
\end{eqnarray*}
and eventually
\begin{eqnarray*}
\mathcal{S}_{i,j} : H^{-\frac12}(\partial B) &\rightarrow& H^{\frac12}(\partial B)\\
\vecf &\mapsto& \int_{\partial B} \vecK(\vecx_i,\vecx_j) \vecf(\vecs)\,d\vecs\,.
\end{eqnarray*}
That these operators are continuous operators from $H^{-\frac12}(\partial B)$ into $H^{\frac12}(\partial B)$ is classical. We hereafter are only interested into the estimation of their norms, and more precisely the way they depend on $a$, $\delta$ and $y$ in the limit $a\rightarrow 0$. Notice that since the kernel $\vecG$ is homogeneous of degree -1, one has 
\begin{equation}
\|{\mathcal{T}}^0\|_{\mathcal{L}(\mathcal{H}^{-1/2},\mathcal{H}^{1/2})} =O\left(\frac{1}{a}\right)\,\mbox{ and }\left\|\left({\mathcal{T}}^0\right)^{-1}\right\|_{\mathcal{L}(\mathcal{H}^{1/2},\mathcal{H}^{-1/2})} =O\left(a\right)\,.
\end{equation}
As far as $\mathcal{S}_l$ is concerned, we get that (since $\ds |\vecK_l(\vecx_i,\vecx_i)| = O\left(1\right)$)
\begin{equation}
\|\mathcal{S}_l\|_{\mathcal{L}(\mathcal{H}^{-1/2},\mathcal{H}^{1/2})} =O\left(1\right)\,,
\end{equation}
and similarly
\begin{equation}
\|\mathcal{S}_{i,j}\|_{\mathcal{L}(\mathcal{H}^{-1/2},\mathcal{H}^{1/2})} =O\left(1\right)\,.
\end{equation}
When $a \rightarrow 0$ this enables us to invert (\ref{expansion}) leading to (\ref{inverse}).
\EndProof

\subsection{Self-propulsion}
 \label{Subsec:Self-propulsion}

We now use the fact that the spheres are non-deformable and may only move following a rigid body motion. In other words, the velocity of each point $\vecr$ of the $i-$sphere expresses as a sum of a translation and a rotation as
\begin{equation}
\vecu_i(\vecr)=\vecu_{T_i} + \vecu_{R_i}(\vecr)\,,
\label{mvtsolide}
\end{equation}
where $\vecu_{T_i}$ is constant on $\partial B$ while $\vecu_{R_i}(\vecr)=\omega_i\times a\vecr$ for a suitable angular velocity $\omega_i$ (remember that all quantities are expressed on the unit sphere $\partial B$).
This is of peculiar importance for the computation of the total force and the total torque, which, due to self-propulsion, should vanish which implies

\begin{equation}
\label{SommeDesForces}
\begin{array}{ll}
\displaystyle \sum_{i} \int_{\partial B} \vecf_i = \displaystyle \sum_{i} \int_{\partial B} \left( \mathcal{T}^{-1}_{\vecx} \vecu \right)_i = 0\,.
\end{array}
\end{equation}
Plugging (\ref{mvtsolide}) in (\ref{SommeDesForces}) and using (\ref{inverse}) leads to
\begin{equation}
\label{TermesDiagonaux}
\begin{array}{l}
\displaystyle \sum_i \int_{\partial B}
({{\mathcal{T}}^0})^{-1} \left({\vecu_T}_i + {\vecu_R}_i - \sum_{k=1}^3  \vecK_k(\vecx_i,\vecx_i) \langle ({{\mathcal{T}}^0})^{-1} ({\vecu_T}_i + {\vecu_R}_i), \mbox{Id} \rangle_{\partial B} \right)-\\
 \hspace*{2.5cm} \displaystyle ({{\mathcal{T}}^0})^{-1} \left(  \sum_{j\neq i} \vecK(\vecx_i,\vecx_j) \langle ({\bar{\mathcal{T}}^0})^{-1} ({\vecu_T}_j + {\vecu_R}_j), \vecId \rangle_{\partial B}\right) = O\left(a^3\right)||\vecu||\,.
\end{array}
\end{equation}

It is well known that both translations and rotations are eigenfunctions of the Dirichlet to Neumann map of the three dimensional Stokes operator outside a sphere. Namely 
$$
\left({{\mathcal{T}}^0}\right)^{-1} {\vecu_T}_i = \lambda_T{\vecu_T}_i\mbox{ and }  \left({{\mathcal{T}}^0}\right)^{-1}{\vecu_R}_i = \lambda_R {\vecu_R}_i\,.
$$
It is well-known that $\lambda_T=\frac{3\mu a}{2}$ leading in particular to the celebrated Stokes formula 
$$
\int_{\partial B}
\left({{\mathcal{T}}^0}\right)^{-1} \vecu_{T_i} \mbox{d}\vecs= 6\pi\mu a \,\vecu_{T_i}
$$
while $\lambda_R = 3\mu a$. We also remark that due to 
$
\int_{\partial B}  \vecu_{R_i} \mbox{d}\vecs = 0\,$, we have
$$
\int_{\partial B}
\left({{\mathcal{T}}^0}\right)^{-1} \vecu_{R_i} \mbox{d}\vecs = 0\,.
$$

We therefore obtain
\begin{equation}
\begin{array}{l}
\displaystyle 6\pi\mu a \sum_i \left({\vecu_T}_i - 6\pi\mu a\sum_{k=1}^3  \vecK_k(\vecx_i,\vecx_i)   {\vecu_T}_i  - 6\pi\mu a \sum_{j\neq i} \vecK(\vecx_i,\vecx_j) {\vecu_T}_j \right) = O\left(a^3\right)||\vecu||\,.
\end{array}
\label{TotalForce}
\end{equation}

We now compute the torque with respect to the center $\vecx_1$ of the first ball $B_1$. Self-propulsion of the swimmer implies that this torque vanishes:
\begin{equation}
\label{TotalTorque}
0 = \int_{\partial B} a\vecr \times \vecf_1(\vecr)+ \int_{\partial B} (\vecx_2-\vecx_1 + a\vecr) \times \vecf_2(\vecr)+\int_{\partial B} (\vecx_3-\vecx_1 + a\vecr) \times \vecf_3(\vecr)  
= \vecI_1 + \vecI_2 + \vecI_3\,,
\end{equation}
where the quantities $\vecI_1$, $\vecI_2$ and $\vecI_3$ are respectively given below.
Calling $\vece_\theta=\left(\begin{array}{c}\cos \theta\\\sin \theta\\ 0\end{array}\right)$ the direction of the swimmer
\begin{eqnarray*}
\vecI_1 &=& \int_{\partial B} (\vecx_1-\vecx_2+a\vecr) \times \vecf_1(\vecr)=  \int_{\partial B} (\xi_1\vece_\theta+a\vecr) \times ({\mathcal{T}}_{\vecx})^{-1} {\vecu_1}\\
&=& \int_{\partial B} \left( -\xi_1 \vece_\theta + a\vecr \right) \times   ({{\mathcal{T}}^0})^{-1}\left({\vecu_T}_1 + {\vecu_R}_1 - 6\pi \mu a \sum_{l=1}^3  \vecK_l(\vecx_1,\vecx_1) {\vecu_T}_1\right.\\
&&\hspace*{5cm}\left. - 6\pi \mu a \sum_{j\neq 2} \vecK(\vecx_1,\vecx_j) {\vecu_T}_j +O\left(a^2\right)||\vecu||\right) \\
&=&-6\pi \mu a  \xi_1 \vece_\theta \times   \left({\vecu_T}_1 - 6\pi \mu a \sum_{l=1}^3  \vecK_l(\vecx_1,\vecx_1) {\vecu_T}_1 - 6\pi \mu a \sum_{j\neq 1} \vecK(\vecx_1,\vecx_j) {\vecu_T}_j\right)+ O\left(a^3\right)||\vecu||\,.\\
\end{eqnarray*}
Similarly, we get,
\begin{eqnarray*}
\vecI_2 &=& a\int_{\partial B} \vecr \times \vecf_2(\vecr)=  a\int_{\partial B} \vecr \times ({\mathcal{T}}_{\vecx})^{-1} {\vecu_2}\\
&=& a\int_{\partial B} \vecr \times   ({{\mathcal{T}}^0})^{-1} \left({\vecu_T}_2 + {\vecu_R}_2 - 6\pi \mu a \sum_{l=1}^3  \vecK_l(\vecx_2,\vecx_2) {\vecu_T}_2 \right.\\
&&\hspace*{5cm}\left.- 6\pi \mu a \sum_{j\neq 2} \vecK(\vecx_2,\vecx_j) {\vecu_T}_j +O\left(a^2\right)||\vecu||\right) \\
&=&a\int_{\partial B} \vecr \times   ({{\mathcal{T}}^0})^{-1} \left({\vecu_R}_1\right) +O\left(a^4\right)||\vecu||\,.
\end{eqnarray*}
But since $\left(\mathcal{T}^0\right)^{-1} {\vecu_R}_1 = \lambda_R {\vecu_R}_1 = \lambda_R \omega_1\times a\vecr$, we have
\begin{eqnarray*}
a\int_{\partial B} \vecr \times   ({{\mathcal{T}}^0})^{-1} \left({\vecu_R}_2\right)  &=& a^2\lambda_R \int_{\partial B} \vecr \times(\omega_1\times \vecr)\,d\vecr\\
&=&\frac{8\pi}{3} \mu a^3 \omega_1\,.
\end{eqnarray*}
This leads to 
\begin{eqnarray*}
\vecI_2 &=& \frac{8\pi}{3} \mu a^3 \omega_1+O\left(a^4\right)||\vecu||\,.
\end{eqnarray*}
Correspondingly, 
\begin{eqnarray*}
\vecI_3 &=& \int_{\partial B} (\vecx_3-\vecx_2+a\vecr) \times \vecf_3(\vecr)\\
&=&6\pi \mu a  \xi_2 \vece_\theta \times   \left({\vecu_T}_3- 6\pi \mu a \sum_{l=1}^3  \vecK_l(\vecx_3,\vecx_3) {\vecu_T}_3  - 6\pi \mu a \sum_{j\neq 3} \vecK(\vecx_3,\vecx_j) {\vecu_T}_j\right)+ O\left(a^3\right)||\vecu||\,.\\
\end{eqnarray*}

Denoting by $\vecA$ the matrix
\begin{equation}
\vecA=\left(
\begin{array}{ccc}
\vecA_{11} & \vecA_{12}  & \vecA_{13}\\
\vecA_{21} & \vecA_{22}  & \vecA_{23}\\
\vecA_{31} & \vecA_{32}  & \vecA_{33}
\end{array}
\right)
\label{defA1}
\end{equation}
where for $i=1,2,3$ 
\begin{equation}
\vecA_{ii}=\vecId - 6\pi\mu a \sum_{l=1}^3 \vecK_l(\vecx_i,\vecx_i)
\label{defA2}
\end{equation}
and for $i,j=1,2,3$ with $i\neq j$
\begin{equation}
\vecA_{ij}=-6\pi\mu a \,\vecK(\vecx_i,\vecx_j)
\label{defA3}
\end{equation}
and $\vecS$ the matrix
$$
\vecS
=
\left(
\begin{array}{ccc}
\vecId & \vecId & \vecId \\
 -\xi_1 \vece_\theta \times  & 0 & +\xi_2 \vece_\theta \times\\
\end{array}
\right)\,,
$$
we can rewrite the self propulsion assumption (\ref{TotalForce}), (\ref{TotalTorque}) as (notice that angular velocities being involved of higher order disappear)
\begin{equation}
\vecS \vecA \left(
\begin{array}{c} \vecu_{T_1}\\ \vecu_{T_2} \\  \vecu_{T_3}
\end{array}\right) = O\left(a^2\right)||\vecu||.
\label{TotalForceTorque}
\end{equation}

We end up by expressing $\vecu_{T_1}, \vecu_{T_2}, \vecu_{T_3}$ and $\omega_1,\omega_2,\omega_3$ in terms of 
$\dot{x},\dot{y},\dot{\theta},\dot{\xi_1}$ and $\dot{\xi_2}$. But, since $\vecu_{T_i}$ is the velocity of the center of the ball $B_i$, one has
$$
\vecu_{T_1} = \left(\begin{array}{c} \dot{x} - \dot{\xi}_1\cos(\theta) +\dot{\theta}\xi_1\sin(\theta) \\ \dot{y} - \dot{\xi}_1\sin(\theta)-\dot{\theta}\xi_1\cos(\theta)\\0 \end{array}\right)\,,\,\,
\vecu_{T_2} = \left(\begin{array}{c} \dot{x} \\ \dot{y} \\0 \end{array}\right)\,,\,\,
$$
and
$$
\vecu_{T_3} = \left(\begin{array}{c} \dot{x} + \dot{\xi}_2 \cos(\theta) -\dot{\theta}\xi_2 \sin(\theta)\\ \dot{y} + \dot{\xi}_2 \sin(\theta)  +\dot{\theta}\xi_2 \cos(\theta)\\0 \end{array}\right)\,.
$$
Similarly
$$
\omega_1=\omega_2=\omega_3=\left(\begin{array}{c} 0\\0\\ \dot{\theta} \end{array}\right)\,.
$$

We rewrite these formulas as
\begin{equation}
\left(\begin{array}{c} \vecu_{T_1}\\ \vecu_{T_2} \\  \vecu_{T_3}
\end{array}\right) = \vecT \left(\begin{array}{c} \dot{x}\\ \dot{y} \\  \dot{\theta}
\end{array}\right) +\vecU \dot{\vecxi}
\label{uti}
\end{equation}
with
$$
\vecT= \left(
\begin{array}{ccc}
\vecId & -\xi_1 \vece_\theta^\perp \\
\vecId & 0 \\
\vecId & \xi_2 \vece_\theta^\perp
\end{array}
\right)\,,
$$
where $\vece_\theta^\perp=\left(\begin{array}{c}-\sin \theta \\ \cos \theta\\ 0\end{array}\right)$
and
$$
\vecU= \left(
\begin{array}{cc}
-\vece_\theta & 0  \\
0 & 0\\
0 & \vece_\theta\\
\end{array}
\right).
$$
Plugging \eqref{uti} into \eqref{TotalForceTorque} leads to the motion equation

\begin{equation}
 \left(\vecS \vecA+\vecR\right)\left( \vecT\left(
 \begin{array}{ll}
\dot{x} \\ \dot{y} \\ \dot{\theta} 
 \end{array}
  \right)+ \vecU \dot{\vecxi}\right) 
 = 0
\label{EquationMouvement1}
 \end{equation}
where the residual matrices have a norm which is estimated as
$$
||\vecR|| =O\left(a^2\right)\,.
$$

\subsection{Dimension of Lie algebra under the small spheres hypothesis}
\label{Subsec:Expansion}
Rewriting from (\ref{defA1}), (\ref{defA2}) and (\ref{defA3}) $\vecA=\vecId+a\vecA_1$, we can expand in power series of $a$ the solution
of (\ref{EquationMouvement1}). This enables us to write an expansion (still in $a$) of the two vectorfields $\vecF_1$ and $\vecF_2$.
To this end, we have used the software MAPLE to symbolically compute those expressions and the Lie brackets $[\vecF_1,\vecF_2]$, $[\vecF_1,[\vecF_1,\vecF_2]]$,
and $[\vecF_1,[\vecF_1,\vecF_2]]$. Writing the vectorfields in components as

\begin{equation}
\label{F1F2}
\vecF_1(\xi_1,\xi_2,y,\theta) := 
\left(
\begin{array}{c}
1 \\ 0  \\ \vecF_1^3 + O\left(a^2\right) \\ \vecF_1^4 +O\left(a^2\right) \\ \vecF_1^5 + O\left(a^2\right)
\end{array}
\right)\, , \vecF_2(\xi_1,\xi_2,y,\theta) := 
\left(
\begin{array}{c}
0 \\ 1  \\ \vecF_2^3 + O\left(a^2\right) \\ \vecF_2^4 +O\left(a^2\right) \\ \vecF_2^5 + O\left(a^2\right)
\end{array}
\right)\, , 
\end{equation}
we find, after having furthermore expanded the abovementioned components in power series of $\ds \frac{1}{y}$,   

\begin{eqnarray*}
\vecF_1^3 &=&  \frac{1}{3} \cos(\theta)+\frac{a}{6}\cos(\theta)K_1^3(\xi_1,\xi_2,\theta)+\frac{3a}{16y^2}\left(\sin(\theta)\cos(\theta)\left(\xi_2+2\xi_1\right)\right)\\
&&+\frac{a}{384y^3}\left(\cos(\theta)K_2^3(\xi_1,\xi_2,\theta)\right)+\frac{a}{512y^4}\left(\sin(\theta)\cos(\theta)K_3^3(\xi_1,\xi_2,\theta)\right)+O\left(\frac{a}{y^5}\right),\\
\vecF_2^3 &=& - \frac{1}{3} \cos(\theta)-\frac{a}{6}\cos(\theta)K_1^3(\xi_2,\xi_1,-\theta)+\frac{3a}{16y^2}\left(\sin(\theta)\cos(\theta)\left(2\xi_2+\xi_1\right)\right)\\
&&-\frac{a}{384y^3}\left(\cos(\theta)K_2^3(\xi_2,\xi_1,-\theta)\right)+\frac{a}{512y^4}\left(\sin(\theta)\cos(\theta)K_3^3(\xi_2,\xi_1,-\theta)\right)+O\left(\frac{a}{y^5}\right)\\
\vecF_1^4 &=& \frac{1}{3} \sin(\theta)+\frac{a}{6}\sin(\theta)K_1^4(\xi_1,\xi_2,\theta)-\frac{3a}{32y^2}K_2^4(\xi_1,\xi_2,\theta)\\
&&+\frac{a}{192y^3}\sin(\theta)K_3^4(\xi_1,\xi_2,\theta)-\frac{a}{y^4} K_4^4(\xi_1,\xi_2,\theta)+O\left(\frac{a}{y^5}\right)\\
\vecF_2^4 &=&  -\frac{1}{3} \sin(\theta)-\frac{a}{6}\sin(\theta)K_1^4(\xi_2,\xi_1,-\theta)-\frac{3a}{32y^2}K_2^4(\xi_2,\xi_1,\theta)\\
&&-\frac{a}{192y^3}\sin(\theta)K_3^4(\xi_2,\xi_1,-\theta)-\frac{a}{y^4}K_4^4(\xi_2,\xi_1,-\theta)+O\left(\frac{a}{y^5}\right)\\
\vecF_1^5 &=&\frac{3a}{64y^3}\sin(\theta)\cos(\theta)K_1^5(\xi_1,\xi_2,\theta)-\frac{9a}{512y^4}\cos(\theta)K_2^5(\xi_1,\xi_2,\theta)+O\left(\frac{a}{y^5}\right)\\ 
\vecF_2^5 &=& \frac{3a}{64y^3}\sin(\theta)\cos(\theta)K_1^5(\xi_2,\xi_1,-\theta)+\frac{9a}{512y^4}\cos(\theta)K_2^5(\xi_2,\xi_1,-\theta)+O\left(\frac{a}{y^5}\right)\,.
\end{eqnarray*} 

In those expressions the remaining functions are respectively given by

\begin{eqnarray*}
K_1^3(\vecxi,\theta)&=&\frac{\left(\xi_2^2\xi_1^2-\xi_2^3\xi_1-\xi_2^4+2\xi_1^3\xi_2+2\xi_1^4\right)}{\left(\xi_2^2+\xi_1\xi_2+\xi_1^2\right)\xi_1\xi_2\left(\xi_1+\xi_2\right)}\\
K_2^3(\vecxi,\theta) &=& -210\xi_1^2\cos(\theta)^2+12\cos(\theta)^4\xi_1^2+184\xi_1^2+24\cos(\theta)^2\xi_1\xi_2\\
&& -32\xi_1\xi_2-6\cos(\theta)^4\xi_1\xi_2-92\xi_2^2+105\xi_2^2\cos(\theta)^2-6\cos(\theta)^4\xi_2^2\\
K_3^3(\vecxi,\theta) &=&\frac{1}{\left(\xi_2^2+\xi_1\xi_2+\xi_1^2\right)}  \Big(12\cos(\theta)^4\xi_2^5+24\xi_1^5\cos(\theta)^4-168\xi_2^5\cos(\theta)^2-336\xi_1^5\cos(\theta)^2\\
& & +112\xi_2^5+72\xi_1\xi_2^4-176\xi_1^2\xi_2^3-136\xi_1^3\xi_2^2+224\xi_1^5-156\xi_1^3\xi_2^2\cos(\theta)^2\\
& &-24\xi_1^2\xi_2^3\cos(\theta)^2-240\xi_1^4\xi_2\cos(\theta)^2+48\xi_1^4\xi_2-156\xi_1\xi_2^4\cos(\theta)^2\\
& & -24\cos(\theta)^4\xi_1^2\xi_2^3+9\cos(\theta)^4\xi_1\xi_2^4-21\xi_1^3\cos(\theta)^4\xi_2^2\Big)\\
K_1^4(\vecxi,\theta) &=&\frac{\left(\xi_2^2\xi_1^2-\xi_2^3\xi_1-\xi_2^4+2\xi_1^3\xi_2+2\xi_1^4\right)}{\left(\xi_2^2+\xi_1\xi_2+\xi_1^2\right)\xi_1\xi_2\left(\xi_1+\xi_2\right)}\\
K_2^4(\vecxi,\theta)&=&6\cos(\theta)^2\xi_1+3\cos(\theta)^2\xi_2-4\xi_1-2\xi_2
\end{eqnarray*}
\begin{eqnarray*}
K_3^4(\vecxi,\theta) &=&-132\xi_1^2\cos(\theta)^2+6\cos(\theta)^4\xi_1^2+56\xi_1^2+12\cos(\theta)^2\xi_1\xi_2\\
 && -16\xi_1\xi_2-3\cos(\theta)^4\xi_1\xi_2-28\xi_2^2+66\xi_2^2\cos(\theta)^2-3\cos(\theta)^4\xi_2^2\\
K_4^4(\vecxi,\theta) &=& \frac{1}{\left(512\xi_2^2+512\xi_1\xi_2+512\xi_1^2\right)}  \Big(-210\cos(\theta)^4\xi_2^5-420\xi_1^5\cos(\theta)^4+232\xi_2^5\cos(\theta)^2\\
&&+24\cos(\theta)^6\xi_1^5-64\xi_2^5+12\cos(\theta)^6\xi_2^5-96\xi_1\xi_2^4-64\xi_1^2\xi_2^3\\
& &-128\xi_1^5+104\xi_1^3\xi_2^2\cos(\theta)^2-56\xi_1^2\xi_2^3\cos(\theta)^2\\
& &-96\xi_1^4\xi_2+216\xi_1\xi_2^4\cos(\theta)^2-66\cos(\theta)^4\xi_1^2\xi_2^3\\
& &-318\xi_1^4\cos(\theta)^4\xi_2-240\xi_1^3\cos(\theta)^4\xi_2^2-24\cos(\theta)^6\xi_1^2\xi_2^3\\
& &-21\cos(\theta)^6\xi_1^3\xi_2^2+464\xi_1^5\cos(\theta)^2-128\xi_1^3\xi_2^2\\
&&+264\xi_1^4\xi_2\cos(\theta)^2-204\cos(\theta)^4\xi_1\xi_2^4+9\cos(\theta)^6\xi_1\xi_2^4\Big),
\end{eqnarray*}

\begin{eqnarray*}
K_1^5(\vecxi,\theta) &=& -8\xi_1-4\xi_2+2\cos(\theta)^2\xi_1+\cos(\theta)^2\xi_2\\
K_2^5(\vecxi,\theta)&=&\frac{1}{\left(\xi_2^2+\xi_1\xi_2+\xi_1^2\right)} \Big(20\cos(\theta)^2\xi_2^4-40\cos(\theta)^2\xi_1^4-4\cos(\theta)^4\xi_2^4+8\xi_1^4\cos(\theta)^4\\
& & -40\xi_2^3\xi_1-8\xi_2^2\xi_1^2+32\xi_1^3\xi_2-16\xi_2^4+32\xi_2^3\cos(\theta)^2\xi_1-7\cos(\theta)^4\xi_2^3\xi_1\\
& &+\cos(\theta)^4\xi_1^2\xi_2^2-40\cos(\theta)^2\xi_1^3\xi_2+8\xi_1^3\cos(\theta)^4\xi_2+32\xi_1^4-8\xi_2^2\xi_1^2\cos(\theta)^2\Big)\,.
\end{eqnarray*}

As one can see, the use of a software for symbolic computation seems unavoidable. Subsequently, we get the expansion of the Lie bracket $[\vecF_1,\vecF_2](\xi_1,\xi_2,y,\theta)$ 
\begin{equation}
\label{ExpressionCrochetLie}
[\vecF_1,\vecF_2](\vecxi,y,\theta) := 
\left(
\begin{array}{cccc}
0 \\ 0 \\ \left[\vecF_1, \vecF_2\right]_3 + O\left(a^2\right) \\ \left[\vecF_1, \vecF_2\right]_4 +O\left(a^2\right) \\ \left[\vecF_1, \vecF_2\right]_5 + O\left(a^2\right)
\end{array}
\right)\, ,
\end{equation}
where the components are given by the following expressions
\begin{eqnarray*}
[\vecF_1, \vecF_2]_3 &=&  -\frac{a}{3}\cos(\theta)\frac{(\xi_1^4+2\xi_1^3\xi_2+\xi_2^2\xi_1^2+2\xi_2^3\xi_1+\xi_2^4)}{(\xi_1+\xi_2)^2\xi_2^2\xi_1^2}\\
&& -\frac{27\,a}{512\,y^4}a\cos(\theta)\sin(\theta)\frac{\xi_1\xi_2\left(\cos(\theta)^4-4\cos(\theta)^2+8\right)\left(-\xi_2^2+\xi_1^2\right)}{\left(\xi_2^2+\xi_1\xi_2+\xi_1^2\right)}\\
&&+ O\left(\frac{a}{y^5}\right)\,,\\
\left[ \vecF_1, \vecF_2\right]_4 &=& -\frac{a}{3}\sin(\theta)\frac{\left(\xi_1^4+2\xi_1^3\xi_2+\xi_2^2\xi_1^2+2\xi_2^3\xi_1+\xi_2^4\right)}{\left(\xi_1+\xi_2\right)^2\xi_2^2\xi_1^2}\\
&&+\frac{27\,a}{512\,y^4}\cos(\theta)^2\frac{\xi_1\xi_2\left(\cos(\theta)^4-4\cos(\theta)^2+8\right)\left(-\xi_2^2+\xi_1^2\right)}{(\xi_2^2+\xi_1\xi_2+\xi_1^2)}\\
&& + O\left(\frac{a}{y^5}\right)\,,
\end{eqnarray*}
\begin{eqnarray*}
\left[\vecF_1, \vecF_2\right]_5 &=&\frac{81\,a}{512\,y^4}\cos(\theta)\frac{\xi_1\xi_2\left(\xi_1+\xi_2\right)\left(\cos(\theta)^4-4\cos(\theta)^2+8\right)}{(\xi_2^2+\xi_1\xi_2+\xi_1^2)} + O\left(\frac{a}{y^5}\right)\,.
\end{eqnarray*}

\noindent Notice that since the two first coordinates of $\vecF_1$ and $\vecF_2$ are constant, the corresponding first coordinates of the Lie bracket vanish.
Similarly, the asymptotic expansion for the second order Lie bracket $\left[\vecF_1,\left[\vecF_1, \vecF_2\right]\right](\vecxi,y,\theta)$ reads

\begin{equation}
\label{ExpressionCrochetLie2}
\left[\vecF_1,\left[\vecF_1, \vecF_2\right]\right](\vecxi,y,\theta) := 
\left(
\begin{array}{c}
0 \\ 0 \\ \left[\vecF_1,\left[\vecF_1, \vecF_2\right]\right]_3 + O\left(a^2\right) \\ \left[\vecF_1,\left[\vecF_1, \vecF_2\right]\right]_4 + O\left(a^2\right) \\ \left[\vecF_1,\left[\vecF_1, \vecF_2\right]\right]_5 + O\left(a^2\right)
\end{array}
\right)\, ,
\end{equation}
where
\begin{eqnarray*}
[\vecF_1,[\vecF_1, \vecF_2]]_3 &=& -\frac{2\,a}{3}\cos(\theta)\frac{\xi_2\left(3\xi_1^2+3\xi_1\xi_2+\xi_2^2\right)}{\xi_1^3(\xi_1+\xi_2)^3}\\
&&+\frac{27\,a}{512\,y^4}\cos(\theta)\sin(\theta)L_3(\vecxi,\theta) + O\left(\frac{a}{y^5}\right) \,,\\
\left[\vecF_1,[\vecF_1, \vecF_2]\right]_4 &=&-\frac{2\,a}{3}\sin(\theta)\frac{\xi_2\left(3\xi_1^2+3\xi_1\xi_2+\xi_2^2\right)}{\xi_1^3\left(\xi_1+\xi_2\right)^3}\\
&&-\frac{27\,a}{512\,y^4}\cos(\theta)^2L_3(\vecxi,\theta)  + O\left(\frac{a}{y^5}\right)\,,\\
\left[\vecF_1,[\vecF_1, \vecF_2]\right]_5 &=&-\frac{81\,a}{512\,y^4}\cos(\theta)L_4(\vecxi,\theta)+ O\left(\frac{a}{y^5}\right)\,.
\end{eqnarray*}
There, $L_3$ and $L_4$ are respectively given by
\begin{eqnarray*}
L_3(\vecxi,\theta) &=& \frac{\xi_2^3\left(8-4\cos(\theta)^2+\cos(\theta)^4\right)\left(2\xi_1^2-\xi_1\xi_2-\xi_2^2\right)}{(\xi_2^2+\xi_1\xi_2+\xi_1^2)^2}\,,\\
L_4(\vecxi,\theta) &=& \frac{\xi_2^3\left(8-4\cos(\theta)^2+\cos(\theta)^4\right)\left(2\xi_1+\xi_2\right)}{(\xi_2^2+\xi_1\xi_2+\xi_1^2)^2}\,.
\end{eqnarray*}

Eventually, the expansion of the vector field $\left[\vecF_2,\left[\vecF_1, \vecF_2\right]\right]$ is given by
\begin{equation}
\label{ExpressionCrochetLie3}
\left[\vecF_2,\left[\vecF_1, \vecF_2\right]\right](\xi_1,\xi_2,y,\theta) := 
\left(
\begin{array}{ccccc}
0 \\ 0 \\ \left[\vecF_2,\left[\vecF_1, \vecF_2\right]\right]_3 + O\left(a^2\right) \\ \left[\vecF_2,\left[\vecF_1, \vecF_2\right]\right]_4 + O\left(a^2\right) \\ \left[\vecF_2,\left[\vecF_1, \vecF_2\right]\right]_5 + O\left(a^2\right)
\end{array}
\right)\, ,
\end{equation}
where
\begin{eqnarray*}
[\vecF_2,[\vecF_1, \vecF_2]]_3 &=& -\frac{2}{3}a\cos(\theta)\frac{\xi_1\left(\xi_1^2+3\xi_1\xi_2+3\xi_2^2\right)}{(\xi_2^3\left(\xi_1+\xi_2\right)^3)}\\
&& -\frac{27\,a}{512\,y^4}a\cos(\theta)\sin(\theta)L_3(\xi_2,\xi_1,-\theta) + O\left(\frac{a}{y^5}\right)\,,\\
\left[\vecF_2,[\vecF_1, \vecF_2]\right]_4 &=&\frac{2\,a}{3}\sin(\theta)\frac{\xi_1\left(\xi_1^2+3\xi_1\xi_2+3\xi_2^2\right)}{(\xi_2^3\left(\xi_1+\xi_2\right)^3}\\
&& -\frac{27\,a}{512\,y^4}\cos(\theta)^2L_3(\xi_2,\xi_1,-\theta) + O\left(\frac{a}{y^5}\right)\,,\\
\left[\vecF_2,[\vecF_1, \vecF_2]\right]_5 &=&-\frac{81\,a}{512\,y^4}\cos(\theta)L_4(\xi_2,\xi_1,-\theta) + O\left(\frac{a}{y^5}\right)\,.
\end{eqnarray*}

\noindent We now can compute an expansion of $\mbox{det}\left(\vecF_1,\vecF_2,[\vecF_1,\vecF_2],\left[\vecF_1,\left[\vecF_1, \vecF_2\right]\right],\left[\vecF_2,\left[\vecF_1, \vecF_2\right]\right]\right)$ which if it does not vanish implies the local controllability of our model swimmer. 
It can be readily checked that we have

\begin{eqnarray}
\mbox{det}\left(\vecF_1,\vecF_2,[\vecF_1,\vecF_2],\left[\vecF_1,\left[\vecF_1, \vecF_2\right]\right],\left[\vecF_2,\left[\vecF_1, \vecF_2\right]\right]\right) &=& \nonumber\\
&&\hspace*{-5cm} = \left| \begin{array}{ccc}
\left[\vecF_1, \vecF_2\right]_3 & \left[\vecF_1,\left[\vecF_1, \vecF_2\right]\right]_3 & \left[\vecF_2,\left[\vecF_1, \vecF_2\right]\right]_3\\
\left[\vecF_1, \vecF_2\right]_4 &\left[\vecF_1,\left[\vecF_1, \vecF_2\right]\right]_4 &\left[\vecF_2,\left[\vecF_1, \vecF_2\right]\right]_4\\ 
\left[\vecF_1, \vecF_2\right]_5 &  \left[\vecF_1,\left[\vecF_1, \vecF_2\right]\right]_5&\left[\vecF_2,\left[\vecF_1, \vecF_2\right]\right]_5
\end{array}
\right|\nonumber \\
\label{determinant}
&&\hspace*{-5cm} = \frac{81\,a^3(\xi_1-\xi_2)}{131072\,y^9}\sin\theta (\cos \theta)^2  R(\vecxi,\theta)  + \, O\left(\frac{1}{y^{10}}\right)\,,
\end{eqnarray}

\noindent with,
\begin{eqnarray*}
R(\vecxi,\theta) &=& \frac{\left(6\xi_1^6+27\xi_1^5\xi_2+50\xi_1^4\xi_2^2+55\xi_1^3\xi_2^3+50\xi_1^2\xi_2^4+27\xi_1\xi_2^5+6\xi_2^6\right)}{\left(\xi_1+\xi_2\right)\left(\xi_2^2+\xi_1\xi_2+\xi_1^2\right)^2\xi_1\xi_2} \\
&&\hspace*{3cm}\times  \left(64-64\cos(\theta)^2+32\cos(\theta)^4-8\cos(\theta)^6+\cos(\theta)^8\right) 
\end{eqnarray*}

It is easily seen that $R$ never vanishes. Therefore, the previous determinant has a non-vanishing first coefficient (in $\ds \frac{1}{y^9}$) which does not vanish for $\xi_1\ne \xi_2$ and $\theta \notin \{0, \frac{\pi}{2},\pi, \frac{3\pi}{2}\}$. Since it is an analytic function of $(\xi_1,\xi_2,y,\theta)$ we deduce that it does not vanish except at most on a negligible set.
This is sufficient to conclude that
$$
d_{(\xi_1,\xi_2,y,\theta)}=5
$$
almost everywhere and the local controllability of the Three-sphere swimmer around such points.

\begin{remark}
Quite strikingly, when $\xi_1=\xi_2$ the first term of the expansion vanishes and one has to go one step further. We find in that case
\begin{eqnarray*}
\mbox{det}\left(\vecF_1,\vecF_2,[\vecF_1,\vecF_2],\left[\vecF_1,\left[\vecF_1, \vecF_2\right]\right],\left[\vecF_2,\left[\vecF_1, \vecF_2\right]\right]\right)(\xi,\xi,y,\theta) &=& 
\, T(\xi,y,\theta) \frac{1}{y^{10}} + \, O\left(\frac{1}{y^{11}}\right)\,,
\end{eqnarray*}
where 
$$
T(\xi,y,\theta) = -\frac{945}{524288}a^3\sin(\theta)^2\cos(\theta)^2\xi\left(\cos(\theta)^4+8-4\cos(\theta)^2\right)^2\,.
$$
This coefficient does not vanish unless $\theta \notin \{0, \frac{\pi}{2},\pi, \frac{3\pi}{2}\}$.  
\end{remark}

\vspace*{0.5cm}

\noindent The case $\theta=0$ or $\pi$. We already know from symmetry that when $\theta=\frac{\pi}{2}$ or $\theta=\frac{3\pi}{2}$, one has $d_{(\xi_1,\xi_2,y,\theta)}\leq 3$. Therefore, it remains to understand the case $\theta=0$ (or $\pi$ by symmetry). The preceding computation does not allow us to conclude about the dimension of the Lie algebra at such points. Indeed, the 2 first coefficients of the expansion of the determinant vanish, and it might well be the case at all orders. Nevertheless, in that case, we can expand the subdeterminant
$$
\Delta = \left|
\begin{array}{cc}
\left[\vecF_1, \vecF_2\right]_3 & \left[\vecF_1,\left[\vecF_1, \vecF_2\right]\right]_3\\
\left[\vecF_1, \vecF_2\right]_5 &\left[\vecF_1,\left[\vecF_1, \vecF_2\right]\right]_5
\end{array}
\right|
$$
in order to obtain informations. Indeed, one gets

\begin{eqnarray*}
\Delta &=& \frac{45}{512\,y^4}a^2\left(\xi_1-\xi_2\right)\xi_1R'(\vecxi)+ O\left(\frac{1}{y^5}\right)\,,
\end{eqnarray*}
\noindent with,

\begin{eqnarray*}
R'(\vecxi) &=& \frac{1}{\left(\xi_1+\xi_2\right)^2\xi_2^2\left(\xi_2^2+\xi_1^2+\xi_1\xi_2\right)^2}  \Big(2\xi_2^5+11\xi_1\xi_2^4+16\xi_1^2\xi_2^3+19\xi_1^3\xi_2^2+12\xi_1^4\xi_2+3\xi_1^5\Big)\,.
\end{eqnarray*}

As the direct consequence, we get that the dimension of the Lie algebra, $d_{\left(\vecxi,y,0\right)}\geq 4$, for almost every $(\vecxi,y)\in (\mathbf{R}^+)^3$. 

This finishes the proof of Lemma \ref{lemmalie} and thus of Theorem \ref{Thm2}.
\EndProof

\begin{remark}
As usual, it is possible to pass from local to global controllability on each of the connected components where the determinant given by (\ref{determinant}) does not vanish. More precisely, let $\mathcal{A} := \left\{\left(\vecxi,\vecp\right) \,s.\,\,t. \, d_{\left(\vecxi,\vecp\right)} \leq 4 \right\} $, we define by $\mathcal{S}_{\left(\vecxi,\vecp\right)}$ the connected component of the subset $\mathcal{S}\setminus \mathcal{A}$ which contains $\left(\vecxi,\vecp\right)$. Applying Chow's Theorem \ref{ChowGlobal} on $\mathcal{S}_{\left(\vecxi,\vecp\right)}$, gives that for every initial configuration $(\vecxi^i,\vecp^i)$, any final configuration $(\vecxi^f,\vecp^f)$ in $\mathcal{S}_{\left(\vecxi,\vecp\right)}$, and any final time $T>0$, there exists a stroke $\vecxi\in \mathcal{W}^{1,\infty}([0,T])$, satisfying
$\vecxi(0)=\vecxi^i$ and $\vecxi(T)=\vecxi^f$ and such that the self-propelled swimmer starting in position $\vecp^i$ with the shape $\vecxi^i$ at time $t=0$, ends at position $\vecp^f$ and shape $\vecxi^f$ at time $t=T$ by changing its shape along $\vecxi(t)$ and staying in $\mathcal{S}_{\left(\vecxi,\vecp\right)}$ for all time $t\in [0,T]$.
In other words, $\mathcal{S}_{\left(\vecxi,\vecp\right)}$ is exactly equal to the orbit of the point $\left(\vecxi,\vecp\right)$.
\end{remark}

\newpage
\section{Conclusion}
\label{Conclusion}

The aim of the present paper was to examine how the controllability of low Reynolds number artificial swimmers is affected  by the presence of a plane boundary on the fluid. The systems are those classically studied in the literature (see \cite{AlougesDeSimone10} for instance) but are usually not confined. This is the first in-depth control study of how the presence of the plane wall affects the reachable set of a peculiar micro-swimmer.

Firstly, the Theorem \ref{Thm1} shows that the controllability on the whole space implies the controllability in the half space. Although the proof is applied on the Four-sphere swimmer, it is based on general arguments which can be appropriate for any finite dimensional linear control systems. 

Secondly, the Theorem \ref{Thm2} deals with the controllability of the Three-sphere swimmer in the presence of the plane wall. We prove that, at least for this example, the hydrodynamics perturbation due to the wall surprisingly makes the swimmer more controllable. This result is not in contradiction with the several scientific studies which show that the wall seems to attract the swimmer (see \cite{Rothschild63}, \cite{WinetBernstein84}, \cite{SmithBlake10}, \cite{GaffneyGadelha11}, \cite{BerkeH.-C.-Berg08}). Although, the Theorem \ref{Thm2} leads to the fact that the wall contributes to increase the swimmer's reachable set, we can conjecture that some of them are easier to reach than others. 

The quantitative approach to this question together than the complete understanding of the situation in view of controllability of the underlying systems is far beyond reach and thus still under progress as is,  
in another direction, the consideration of more complex situations like, e.g. rough or non planar wall. This is the purpose of ongoing work. 

\textbf{Acknowledgements} :
This work has been supported by Direction Générale de l'Armement (DGA).
\newpage

\bibliographystyle{plain}


\begin{thebibliography}{}

\end{thebibliography}


\begin{thebibliography}{10}

\bibitem{AlougesDeSimone10}
F.~Alouges, A.~DeSimone, L.~Heltai, A.~Lefebvre and B.~Merlet.
\newblock Optimally swimming stokesian robots.
\newblock {\em Preprint arXiv :1007.4920v1 [math.OC]}, 2010.

\bibitem{AlougesDeSimone08}
F~Alouges, A.~DeSimone and A.~Lefebvre.
\newblock Optimal strokes for low Reynolds number swimmers : an example.
\newblock {\em Journal of Nonlinear Science}, 18:277--302, 2008.

\bibitem{BerkeAllison08}
A.~P. Berke and P.~Allison.
\newblock Hydrodynamic attraction of swimming microorganisms by surfaces.
\newblock {\em Physical Review Letters}, 101(038102), 2008.

\bibitem{BerkeH.-C.-Berg08}
A.~P. Berke, L.~Turner, H.~C. Berg and E.~Lauga.
\newblock Hydrodynamic attraction of swimming microorganisms by surfaces.
\newblock {\em Physical Review Letters}, 2008.

\bibitem{Blake71}
J.~R. Blake.
\newblock A note on the image system for a Stokeslet in a no-slip boundary.
\newblock {\em Proc. Gamb. Phil. Soc.}, 70:303, 1971.

\bibitem{Chambrion11}
T. Chambrion and A. Munnier.
\newblock Generic controllability of 3D swimmers in a perfect fluid.
\newblock {\em Preprint arXiv:1103.5163,}  2011.

\bibitem{Coron56}
J.~M. Coron.
\newblock {\em Control and Nonlinearity}.
\newblock Mathematical Surveys and Monographs, Vol. 136, American Mathematical Society, Providence, 2007.

\bibitem{Dal-MasoDeSimone10}
G.~Dal~Maso, A.~DeSimone and M.~Morandotti.
\newblock An existence and uniqueness result for the selfpropelled
motion of micro-swimmers.  
\newblock {\em SIAM J. Math. Anal.}, 43:1345-1368, 2011.

\bibitem{GaffneyGadelha11}
E.~A. Gaffney, H.~Gad{\^e}lha, D.~J. Smith, J.~R. Blake and J.~C.
  Kirkman-Brown.
\newblock Mammalian sperm motility: observation and theory.
\newblock {\em J. Fluid Mech.}, 2011.

\bibitem{GolestanianAjdari08}
R.~Golestanian and A.~Ajdari.
\newblock Analytic results for the three-sphere swimmer at low Reynolds.
\newblock {\em Physical Review E}, 77:036308, 2008.

\bibitem{Jurdjevic97}
V.~Jurdjevic.
\newblock {\em Geometric Control Theory}.
\newblock Cambridge University Press., 1997.

\bibitem{LeshanskyKenneth08}
A.~M. Leshansky and O.~Kenneth.
\newblock Surface tank treading: Propulsion of Purcell's toroidal swimmer.
\newblock {\em Physics of Fluids}, 20(063104), 2008.

\bibitem{LoheacMunnier12}
J.~Loh{\'e}ac, and A.~Munnier.
\newblock Controllability of 3D low Reynolds swimmers.
\newblock {\em Preprint arXiv:1202.5923}, 2012.

\bibitem{LoheacScheid11}
J.~Loh{\'e}ac, J.~F. Scheid and M.~Tucsnak.
\newblock Controllability and time optimal control for low Reynolds numbers
  swimmers.
\newblock {\em Preprint Hal 00635981}, 2011.

\bibitem{Munnier09}
A.~Munnier.
\newblock Locomotion of deformable bodies in an ideal fluid: Newtonian versus
  lagrangian formalisms.
\newblock {\em J. Nonlinear Sci.}, 2009.

\bibitem{NajafiGolestanian04}
A.~Najafi and R.~Golestanian.
\newblock Simple swimmer at low Reynolds number: Three linked spheres.
\newblock {\em Physical Review E}, 69(6):062901, 2004.

\bibitem{NajafiZargar10}
A.~Najafi and R.~Zargar.
\newblock Two-sphere low Reynolds-number propeller.
\newblock {\em Physical Review E}, 81, 2010.

\bibitem{OrMurray09}
Y.~Or and M.~Murray.
\newblock Dynamics and stability of a class of low Reynolds number swimmers
  near a wall.
\newblock {\em Physical Review E}, 79:045302(R), 2009.

\bibitem{Purcell77}
E.~M. Purcell.
\newblock Life at low Reynolds number.
\newblock {\em American Journal of Physics}, 45:3--11, 1977.

\bibitem{Rothschild63}
L.~Rothschild.
\newblock Non-random distribution of bull spermatozoa in a drop of sperm
  suspension.
\newblock {\em Nature}, 1963.

\bibitem{Sauvage01}
J. P. Sauvage.
\newblock {\em Molecular Machines and Motors}.
\newblock Springer, 2001.

\bibitem{ShumGaffney10}
H.~Shum, E. A. Gaffney and D. J.~Smith.
\newblock Modeling bacterial behavior close to a no-slip plane boundary : the
  influence of bacterial geometry.
\newblock {\em Proceeding of the Royal Society A}, 466(1725-1748), 2010.

\bibitem{SmithBlake10}
D.~J. Smith and J.~R. Blake.
\newblock Surface accumulation of spermatozoa: a fluid dynamic phenomenon.
\newblock {\em Preprint arXiv: 1007.2153v1}, 2010.

\bibitem{SmithGaffney09}
D.~J.~Smith, E.~A.~Gaffney, J.~R.~Blake and J.~C.~Kirkman-Brown.
\newblock Human sperm accumulation near surfaces : a simulation study.
\newblock {\em J. Fluid Mech.}, 621:289-320, 2009.

\bibitem{Taylor51}
G.~Taylor.
\newblock Analysis of the swimming of microscopic organisms.
\newblock {\em Proc. R. Soc. Lond. A}, 209:447--461, 1951.

\bibitem{WatsonFriend09}
B.~Watson, J.~Friend and L.~Yeo.
\newblock Piezoelectric ultrasonic resonant motor with stator diameter less
  than 250$\mu$m: the proteus motor.
\newblock {\em J. Micromech. Microeng.}, 19, 2009.

\bibitem{WinetBernstein84}
H.~Winet, G.~S.~Bernstein and J.~Head.
\newblock observation on the response of human spermatozoa to gravity,
  boundaries and fluid shear.
\newblock {\em Reproduction}, 70, 1984.

\bibitem{ZargarNajafi09}
R.~Zargar, A.~Najafi and M.~Miri.
\newblock Three-sphere low Reynolds number swimmer near a wall.
\newblock {\em Physical Review E}, 80:026308, 2009.

\end{thebibliography}

\end{document}